\documentclass[12 pt]{article}
\usepackage{hyperref,amsmath,amsfonts,graphicx,psfrag,fancybox,qsymbols,indentfirst}
\usepackage[latin1]{inputenc}
\def\d{\displaystyle}

\def \sous#1#2{\mathrel{\mathop{\kern 0pt#1}\limits_{#2}}}
\def \sur#1#2{\mathrel{\mathop{\kern 0pt#1}\limits^{#2}}}

\def \be{\begin{eqnarray*}}
\def \ee{\end{eqnarray*}}
\def \ben{\begin{eqnarray}}
\def \een{\end{eqnarray}}
\def\d{\displaystyle}

\def\AArm{\fam0 }
\def\AAk#1#2{\setbox\AAbo=\hbox{#2}\AAdi=\wd\AAbo\kern#1\AAdi{}}%
\def\AAr#1#2#3{\setbox\AAbo=\hbox{#2}\AAdi=\ht\AAbo\raise#1\AAdi\hbox{#3}
}%

\def \sur#1#2{\mathrel{\mathop{\kern 0pt#1}\limits^{#2}}}
\def\BBc{{\AArm C\AAk{-1.02}{C}\AAr{.9}{I}{\AAFf\char"3F}}}%
\def\BBd{{\AArm I\!D}}%
\def\BBe{{\AArm I\!E}}%
\def\BBh{{\AArm I\!H}}%
\def\BBn{{\AArm I\!N}}%
\def\BBr{{\AArm I\!R}}%
\def\BBone{{\AArm 1\AAk{-.8}{I}I}}%

\def \e{{\rm e}}

\def \d{{\tt d}}

\def\noi{\noindent\underbar}

\def\videbox{\mathbin{\vbox{\hrule\hbox{\vrule height1ex \kern.5em\vrule height1ex}\hrule}}}

\def \build#1#2#3{\mathrel{\mathop{\kern 0pt#1}\limits_{#2}^{#3}}}
\def\suPP#1#2{{\displaystyle\sup _{\scriptstyle #1\atop \scriptstyle #2}}}
\def\proDD#1#2{{\displaystyle\prod _{\scriptstyle #1\atop \scriptstyle #2}}}
\def\proof{\noindent{\bf Proof:}\hskip10pt}
\def\QED{\hfill\vrule height 1.5ex width 1.4ex depth -.1ex \vskip20pt}

\begin{document}
\newtheorem{lem}{Lemma}[section]
\newtheorem{defi}[lem]{Definition}
\newtheorem{theo}[lem]{Theorem}
\newtheorem{cor}[lem]{Corollary}
\newtheorem{prop}[lem]{Proposition}
\newtheorem{nota}[lem]{Notation}
\newtheorem{rem}[lem]{Remark}
\newqsymbol{`P}{\mathbb{P}}
\newqsymbol{`E}{\mathbb{E}}
\newqsymbol{`N}{\mathbb{N}}
\newqsymbol{`C}{\mathbb{C}}
\newqsymbol{`R}{\mathbb{R}}
\newqsymbol{`S}{\mathbb{S}}
\newqsymbol{`o}{\omega}
\newqsymbol{`M}{{\cal M}}
\newqsymbol{`H}{\mathbb{H}}
\newqsymbol{`X}{\mathbb{X}}
\newdimen\AAdi%
\newbox\AAbo%
\font\AAFf=cmex10

\setcounter{page}{1}
\setcounter{section}{0}\date{}

\def \build#1#2#3{\mathrel{\mathop{\kern 0pt#1}\limits_{#2}^{#3}}}
\def\suPP#1#2{{\displaystyle\sup _{\scriptstyle #1\atop \scriptstyle #2}}}
\def\proDD#1#2{{\displaystyle\prod _{\scriptstyle #1\atop \scriptstyle #2}}}
\def\proof{\noindent{\bf Proof:}\hskip10pt}
\def\QED{\hfill\vrule height 1.5ex width 1.4ex depth -.1ex \vskip20pt}
\setcounter{page}{1}
\setcounter{section}{0}
\title{Pathwise asymptotic behavior of 
random determinants in the uniform Gram and Wishart ensembles
}
\large
\author{A. Rouault \thanks
{LAMA, B\^atiment Fermat,
Universit\'e de Versailles F-78035 Versailles. e-mail:
rouault@math.uvsq.fr}
}
\maketitle
\normalsize\rm
\centerline{\today}
\bigskip
\noindent{\bf Summary. }{\small This paper concentrates on asymptotic properties of 
determinants of some random symmetric matrices. If $B_{n,r}$ is a $n\times r$ rectangular matrix 
and $B_{n,r}'$ its transpose, we study $\det (B_{n,r}'B_{n,r})$ when $n,r$ tends to infinity with 
$r/n \rightarrow c\in (0,1)$. The $r$ column vectors of $B_{n,r}$ are chosen independently, with common distribution $\nu_n$. 
The Wishart ensemble corresponds to $\nu_n = {\cal N}(0 , I_n)$, the standard normal distribution. We call 
uniform Gram  ensemble the ensemble corresponding
to $\nu_n = \sigma_n$, the uniform distribution on the unit sphere  $`S_{n-1}$. 
In the Wishart ensemble, a well known Bartlett's theorem  decomposes  the above determinant 
into a product of chi-square variables. The same holds in the uniform Gram ensemble. This allows us to study the {\it process} 
 $\{\frac{1}{n}\log \det\big(B_{n,\lfloor nt\rfloor}'B_{n,\lfloor nt\rfloor}\big) , t \in [0,1]\}$ and its 
  asymptotic behavior as $n\rightarrow \infty$: a.s. convergence, fluctuations, large deviations. 
  We connect the results for marginals (fixed $t$) with those obtained by the spectral method.}
 \vskip 5mm
\noindent
{\bf A.M.S. Classification.}  {\tt 15 A 52, 15 A 15, 60F 10, 60F 17, 62 H 10}
\vskip 3mm

\normalsize\rm
\section{Introduction} 
For $n, r \in \BBn$  such that  $r\leq n$,   
let $M_{n,r}(`R)$ be the set of $n\times r$ matrices with real entries.
A matrix $B\in M_{n,r}(`R)$ consists in 
$r$ column vectors $b_1, \cdots , b_r$ of $`R^n$. We denote by  $B'$ its transpose, so that $B'B \in M_{r,r}$ is symmetric. 
We provide $`R^n$ with the usual Euclidean norm.
 
In 1893,  Hadamard  \cite{Hadamard} proved that
\be
\det (B'B) \leq \Vert b_1\Vert^2 \cdots\Vert b_r\Vert^2
\ee 
with equality if and only if  $b_1, \cdots , b_r$ are orthogonal. That means that the volume (or $r$-content) 
of the parallelotope built from  $b_1, \cdots, b_r$ is maximal when the vectors are orthogonal.
Consequently, the quantity 
\ben
\label{rh}
h(B) =\frac{\det (B'B)}{\Vert b_1\Vert^2 \cdots\Vert b_r\Vert^2}
\,\een
is usually called the Hadamard ratio (cf. \cite{dixon}); in the basis reduction problem
(\cite{akhavi1},\cite{akhavi2},\cite{akhavi3}), 
the quantity $1/\sqrt{h(B)}$ is called  the orthogonality defect.
Some papers (\cite{abbott},\cite{dixon}) are concerned with the 
tightness of the bound  $h(B) \leq 1$ when 
$B$ is random and $n=r$.
Writing  $B_{n,r}$ instead of $B$ to stress on  dimensions,
it is interesting to study
 the asymptotic behavior of the sequence of random variables  $h (B_{n,r})$, in particular 
when $n, r \rightarrow \infty$ with $r/n\rightarrow c\in [0,1]$. 

We consider  independent random vectors $b_i, i = 1, \cdots ,r$ with the same distribution $\nu_n$ in $`R^n$. 
It seems natural to choose  $\nu_n = \sigma_n$, the uniform distribution on the unit sphere $`S_{n-1}$. 
The corresponding ensemble for $B$ is called Uniform Spherical Ensemble in \cite{Donoho1}.
The matrix ensemble for $B'B$ is called the Gram ensemble in \cite{Spince1} 
since $B'B$ is the Gram matrix built from vectors $b_i$'s. To stress on the distribution, we call it uniform Gram ensemble. 

More generally,  if $\nu_n$ is isotropic, (i.e. $\nu_n (\{0\})=0$ and $\nu_n$ invariant by rotation), it is well known that
 $\widetilde b_1 := b_1/ \Vert b_1\Vert$ is $\sigma_n$ distributed and  independent of $\Vert b_1 \Vert$. 
Denoting by $\widetilde B$ the matrix of unitary vectors, we see that $\widetilde B'\widetilde B$ is in the uniform Gram ensemble. 
 It makes possible to study $\det B'B$ in its own, since the decomposition  in independent factors
\ben
\label{tilde}
\det (B'B) = \det (\widetilde B' \widetilde B) \times \prod_{i=1}^r \Vert b_i\Vert^2
\een
 reduces this case to the previous one if the distribution of $\Vert b_1\Vert^2$ is well behaved.
  
The most important example is  the Gaussian one with
 $\nu_n = {\cal N}(0 ; 
I_n)$: all the entries of $B$ are i.i.d. ${\cal N}(0 ; 1)$
  and  $B'B$ is in the Wishart ensemble. Moreover $\Vert b_1\Vert^2$ is $\chi_n^2$ distributed. Our paper is concerned essentially with these two cases.

  We introduce a probability space on which all uniform Gram and Wishart matrices are defined simultaneously.
   It is just the infinite product space
   generated by a double infinite sequence of  i.i.d. ${\cal N}(0;1)$ variables $\{b_{i,j}\}_{i,j =1}^\infty$. 
   Then we take 
   $B_{n,r} = \{ b_{i,j},  i=1 , \cdots r , j=1 , \cdots r\}$ and 
  omitting the dimension index $n$, we set
  $\widetilde b_{i,j} = b_{i,j}/ (\sum_{k=1}^n b_{k,j}^2)^{1/2}$  and   
  $\widetilde B_{n,r} = \{ \widetilde b_{i,j},  i=1 , \cdots r , j=1 , \cdots r\}$.

In Section \ref{dec}, we recall some known results. Using the  classical $QR$ decomposition of  $B$ with $Q\in M_{n,r}$ orthogonal and $R\in M_{r,r}$ uppertriangular 
(\cite{Bhatia}), we get $$\det (B_{n,r}'B_{n,r}) = \prod_{j=1}^r R_{jj}^2\,.$$  
In the Wishart case, the variables $R^2_{jj}, j=1, \cdots, r$ are independent  and  $\chi^2$ distributed 
with respective parameters $n-j+1, j=1, \cdots, r$. This result is known as the celebrated  Bartlett decomposition.
In the Gram case, the corresponding variables $\widetilde R^2_{jj}, j=2, \cdots, r$ are independent  and  
beta distributed with respective parameters $\big(\frac{n-j+1}{2}, \frac{j-1}{2}\big)$. 
Therefore we will consider
 $\left\{ \log \det (B_{n,r}'B_{n,r}), r = 1, \cdots , n\right\}$ and its "tilde" version
as  triangular arrays 
  and  prove pathwise\footnote{We stress that this study is pathwise in the parameter $t = r/n$ 
 and not in the "time" parameter as in Wishart processes defined from Brownian matrices. 
} results for the sequence of processes
$\left\{ \frac{1}{n}\log \det (B_{n, \lfloor nt\rfloor}'B_{n, \lfloor nt\rfloor}), \ t\in [0,1] \right\}$.

In Section \ref{MP}, we present  the spectral approach. It starts from 
\begin{eqnarray*}
\log \det   (B_{n,r}'B_{n,r}) = \sum_{k=1}^r \log \lambda_{n,r}^{(k)} 
 \ , \ 
\end{eqnarray*}
where $\lambda_{n,r}^{(k)}, k=1, \cdots , r$ are the (real) eigenvalues of $B_{n,r}'B_{n,r}$. 
We may take advantage of known results (recalled in Section \ref{MP})
on the convergence of the 
empirical spectral distribution (ESD) to the Mar\v cenko-Pastur distribution 
as $r/n \rightarrow c \in (0,1)$ (\cite{MarchPastur1} for the Wishart ensemble and \cite{Jiang2}, \cite{Spince1} 
for the uniform Gram ensemble). 
This allows in the following sections to  recover results for marginals (only) and in the Wishart case to discover fruitful connection between the two methods. 
 
 In Section \ref{Hr}, we study the Gram ensemble and set
 \be
 {\cal Y}_{n,r} = \widetilde B_{n,r}'\widetilde B_{n,r}\ , \ \  \Upsilon_{n,r} = \log \det {\cal Y}_{n,r}\,.
 \ee 
 We state the a.s. convergence 
 of $\left\{ \frac{1}{n}
   \Upsilon_{n, \lfloor nt\rfloor}, \ t\in [0,1] \right\}$ (Theorem \ref{theocvps}),
 the weak convergence of fluctuations 
 (Theorem \ref{DonskerH}) and
 a large deviation principle  (Theorem \ref{LDPH} and Theorem \ref{margH}).
 
 In Section \ref{Wm}, we study the Wishart ensemble.
Since $`E \Vert b_1\Vert^2 = n$, it is natural to normalize  Wishart matrices and set 
\ben
\label{defwh}
{\cal X}_{n,r} = \frac{1}{n}B'_{n,r}B_{n,r}\ , \ \  \Xi_{n,r} = \log \det {\cal X}_{n,r}\,.
\een
The asymptotic behavior of the process
 $\left\{\frac{1}{n} \Xi_{n, \lfloor nt\rfloor}, \ t\in [0,1] \right\}$ 
 is easily deduced from (\ref{tilde}) and the above results in the uniform Gram case. (Of course, it is also possible to use the Bartlett decomposition).

 Section \ref{Je} is devoted to some remarks about extensions to matrices with entries with Gaussian entries
  in other fields (complex, quaternionic), and even with non Gaussian entries. 
An extension to the Jacobi ensemble will be considered in a forthcoming paper.

The proofs of results are located in Section \ref{sect7} and \ref{sect8}.

 All along the paper we use the function
$\ell(x) = \log \Gamma(x)$, and its derivative $\Psi = \ell'$ which   
is the digamma function. Some useful properties of $\ell$ and $\Psi$ are given in Appendix. 
We use also the following functions:
\be
{\cal J}(u)&=& u\log u -u + 1, \ \hbox{for}\ u > 0\\
\nonumber
{\cal J}(u)&=& 1, \ \hbox{for}\ u=0\\
\nonumber
{\cal J}(u) &=& +\infty, \ \hbox{for} \ u<0\,,
\ee
and for $t \geq 0$,
\be
F(t) = \int_0^t {\cal J}(u)\ du = \frac{t^2}{2}\log t -\frac{3t^2}{4} + t\,.
\ee
Recall also that for $a >0$, the $\chi_a^2$ distribution has density
$$\frac{x^{\frac{n}{2} -1}}
{2^{\frac{n}{2}}\Gamma\left(\frac{n}{2}\right)} e^{-\frac{x}{2}}\ \ (x > 0)$$
and that, for $\alpha > 0, \ \beta >0$ the beta$(\alpha , \beta)$ distribution has density
$$\frac{\Gamma(\alpha + \beta)}{\Gamma(\alpha)\Gamma(\beta)} \!\ x^{\alpha -1} (1-x)^{\beta -1} \ \ (x > 0)$$
Let us end this introduction with some comments.
 Wishart matrices were first 
introduced in multivariate statistical analysis as sample covariance matrices 
: $B \in M_{n,r}$  is a data matrix where $r$ is the number of variates and $n$ is the sample size. 
Then $\det (B'B)$ is Wilks' generalized variance (up to a transformation). 
It is used to build tests on the covariance matrix (\cite{Anderson}, \cite{Muir}). 
In contemporary multivariate analysis, it is common to consider large $r$ and large $n$ (see examples in \cite{john}) 
, although it may seem non standard to increase the number of variables for a given size of the sample. Besides, 
in stochastic geometry it seems interesting to describe the evolution of the $r$-content 
of a random $r$-parallelotope as $r$ increases. In quantum dynamics, the uniform Gram ensemble  is introduced by  De Conck 
et al. in \cite{Spince2} and they called $c$ (the limit of $r/n$) a time-parameter, although they assumed it fixed.

\section{Decompositions}
\label{dec}
This section consists in notation and recalls. The key point is a  
decomposition of  determinants in products and its consequence for random Gaussian matrices.

\subsection{Some linear algebra}
\label{linalg}
Every 
matrix $B \in M_{n, r}$ may be decomposed (see \cite{Bhatia}) in a product $B= QR$ of an uppertriangular matrix
$R\in M_{r,r}(`R)$ and an orthogonal matrix $Q\in M_{n,r}(`R)$. 
If  the vectors $b_i, \ i= 1, \cdots , r$ are linearly independent, the decomposition is unique if we 
force diagonal elements of $R$ to be positive.
 By the Gram-Schmidt method,
we set $c_1 = b_1$ and for $2\leq j\leq r$
\be
c_j = b_j - \sum_{k = 1}^{j-1} \frac{<c_k, b_j>}{\Vert c_k\Vert^2}\ c_k\,,   
\ee
and then build the orthonormal system:
\be
f_j = \frac{c_j}{\Vert c_j\Vert}, \ \ 1\leq j\leq r\,.
\ee
This yields
\ben
\label{pytha}
b_j = <f_j , b_j> f_j + \sum_{k = 1}^{j-1} <f_k, b_j> f_k
\een
Now $Q$ consists in $f_1, \cdots  , f_r$ and  $R$ is given by 
\ben
\label{tjj}
R_{j,j} = \Vert c_j\Vert = <f_j , b_j>
, \ \ 1\leq j\leq r\,,
\een
and for $2\leq j\leq r$ and $k \leq j -1$:  
\ben
\label{6}
 R_{k,j} = <f_k, b_j> = \frac{<c_k , b_j>}
{\Vert c_k\Vert}\,.
\een
From (\ref{pytha}) we deduce
\ben
\label{12}
\Vert b_j\Vert^2 = R_{j,j}^2 + \sum_{k=1}^{j-1} |R_{k,j}|^2\,.
\een
We can write $b_i = \Vert b_i\Vert \widetilde b_i$  with 
$\widetilde b_i \in `S_{n-1}$, so that
$f_1, \cdots , f_r$ depend only upon 
$\widetilde b_1, \cdots, \widetilde b_r$. We have  
\ben
\label{Rbr}
R_{jj}^2 = \Vert b_j\Vert^2 \widetilde R_{jj}^2 \ \ , \ \ \widetilde R_{jj}^2 := <f_j , \widetilde b_j>^2\,.
\een
Since $R$ is upper triangular and $B'B = R'R$ we get easily
\ben
\label{bart}
\det (B'B) =   \prod_{j=1}^r R_{jj}^2\,,
\een
and from (\ref{rh}) and (\ref{Rbr})
\ben
\label{h=theta}
 h (B) = \det (\widetilde B'\widetilde B) = \prod_{j=2}^r \widetilde R_{jj}^2\,.
\een
(It is clear, of course that $h(B) \leq 1$, as Hadamard noticed).
\subsection{Random Gaussian Matrices and Bartlett's decomposition}
\label{bartlettr}
In the sequel, we study  models of random matrices in which all entries are independent and ${\cal N}(0,1)$ distributed. 

We are in the situation of Section \ref{linalg}. It is clear that for every $j=1, \cdots, r$,
$$
\hbox{Span}\ \{b_1 , \cdots, b_r\} = \hbox{Span} \ \{c_1 , \cdots, c_r\}= \hbox{Span} \ \{f_1 , \cdots, f_r\} 
\,,$$
and $R_{j,j}$ is a measurable function of $(b_1 , \cdots , b_j)$ thanks to (\ref{tjj}) and (\ref{6}). 

\begin{prop}
\label{celebrated}
\begin{itemize}
\item[1)] If $1\leq r \leq n$,
the random variables $R_{j,j}^2 , \ j= 1, \cdots , r$ are independent and 
\be
R_{j,j}^2\buildrel{\cal D}\over{=}
\chi^2_{n-(j-1)}\, ,
\ee
where $\buildrel{\cal D}\over{=}$ stands for equality in distribution.
 \item[2)] If $2\leq r \leq n$,
the random variables  $\widetilde R_{j,j}^2, \ j=2, \cdots, r$ are independent and
\be
\widetilde R_{j,j}^2 \buildrel{\cal D}\over{=}   beta\left(\frac{n-j+1}{2}, \frac{j-1}{2}\right)\,.
\ee
\end{itemize}
\end{prop}
The first claim is the celebrated Bartlett decomposition (\cite{Bart}). It is quoted in many books and articles in particular
 \cite{Anderson} pp. 170-172, \cite{Muir} pp. 99 th. 3.2.14, \cite{Kshi2}.  The second claim comes from 1) and 
 (\ref{Rbr}). 
For the sake of completeness, we give here the proof of 1), with the so-called "random orthogonal transformation", 
which may be found in \cite{Kshi1}.
\medskip

\proof
1) Let us fix $2 \leq j\leq r$, and condition upon $(b_1 , \cdots , b_{j-1})$. From equation (\ref{6}) we have 
${\cal R}_j := \big(R_{1,j}, \cdots , R_{(j-1), j}\big)' = Fb_j$ where 
 $F = \big(f_1 , \cdots , f_{j-1}\big) \in M_{n, j-1}$ is a (known) orthogonal matrix. 
The Cochran theorem and (\ref{12}) imply that 
${\cal R}_j$ is ${\cal N}(0; I_{j-1})$ distributed, 
that $\Vert {\cal R}_j\Vert^2 = \sum_{k=1}^{j-1} |R_{k,j}|^2$ is 
$\chi^2_{j-1}$ distributed and that 
$R_{jj}^2 =\Vert b_j\Vert^2 - \Vert {\cal R}_j\Vert^2$ is independent of $\Vert {\cal R}_j\Vert^2$ and 
$\chi^2_{n-(j-1)}$ distributed. 
Since in all the above statements, the conditioning variables $(b_1, \cdots , b_{j-1})$ 
did not appear, these statements are true unconditionally. In particular, $R_{j,j}^2$ is independent of $(b_1 , \cdots , b_{j-1})$. This yields that 
all the variables $R_{j,j}^2 , \ j= 1, \cdots , r$ are independent.

2) From (\ref{Rbr}) and the previous remarks, we see that the variables $\widetilde R_{jj}^2$ are also independent. 
To get the distributions, recall that
 if $U\buildrel{\cal D}\over{=}\chi^2_a$ and $V \buildrel{\cal D}\over{=}\chi^2_b$ are 
independent variables,
 then  
$U/(U+V) \buildrel{\cal D}\over{=} beta(a/2, b/2)$.
\QED

To consider the asymptotic behavior of uniform Gram 
 and Wishart determinants in a dynamic (or pathwise) way, let us give some notation.

Set $ \Upsilon_{n,0} = 0, \  \Upsilon_{n,1} = 0$ 
and for $2\leq r \leq n$
\ben
\label{defXi}
 \Upsilon_{n,r}
  =   \sum_{k=2}^{r} \log \widetilde R_{j,j}^2
 \een
which provides a first triangular array. Besides, from (\ref{defwh}) and (\ref{bart}), for $r = 1, \cdots , n$
\ben
\label{secondarray}
\Xi_{n,r} =  \sum_{j=1}^{r} \log\frac{R_{jj}^2}{n}\,.
\een
provides a second triangular array. Actually, in that case, we have also,  from (\ref{Rbr}) and (\ref{h=theta}):
\ben
\label{HW}
\Xi_{n,r}
=  \Upsilon_{n,r}  + S_{n,r}\,,\ \ r = 1, \cdots , n\,,
\een
where 
\be
S_{n,r} := \sum_{k=1}^r \log \frac{\Vert b_k\Vert^2}{n}\,,\ \ r = 1, \cdots , n\,.
\ee
In this auxiliary triangular array,
the independent variables $\big(\Vert b_k\Vert^2, k= 1, \cdots , n
\big)$ are independent of $\big(\Xi_{n,r}, r= 1, \cdots, n\big)$ and
$\Vert b_k\Vert^2 \buildrel{\cal D}\over{=} \chi^2_(n)$. 
The three processes are denoted by
\be
 \Upsilon_n (t) :=  \Upsilon_{n, \lfloor nt\rfloor} \ , \ \Xi_n (t) = \Xi_{n, \lfloor nt \rfloor} \ , \ S_n (t) = S_{n, \lfloor nt \rfloor}
\ , \ t \in [0,1]\,.
\ee

\section{The spectral method}
\label{MP}
Beside the above "decomposition method" we will use the spectral approach which we describe now.

 Let $\lambda_{n,r}^{(k)}, k=1, \cdots , r$ be the (real) eigenvalues of
${\cal X}_{n,r}$ in the regime
$n,r \rightarrow \infty$ with $r/n \rightarrow c < 1$ fixed.
We set
\be
\mu_{n,r} = \frac{1}{r} \sum_{k=1}^r \delta_{\lambda_{n,r}^{(k)}}
\ee
the empirical spectral distribution (ESD). In particular
\be
\int \log x \!\ d\mu_{n,r}(x) 
= \frac{1}{r} \sum_{k=1}^r \log \lambda_{n,r}^{(k)} = \frac{1}{r} \log \det {\cal X}_{n,r}\,.
\ee

For $c > 0$ and $\sigma > 0$, let $\pi_{\sigma^2}^c$ be the probability distribution on $`R$ defined by
\ben
\label{defpi}
\pi_{\sigma^2}^c (dx)= (1-c^{-1})_{_+} \delta_0 (dx) + \frac{\Big((x - \sigma^2 a(c))(\sigma^2 b(c) - x)\Big)_+^{1/2}}
{2\pi \sigma^2c x}\ dx\,,
\een
where $\delta_0$ is the Dirac mass in $0$,  $x_{_+} = \max (x,0)$ and 
\ben
\label{aetb}
a(c) = (1 - \sqrt{c})^2 \ , \ \ b(c) = (1 + \sqrt{c})^2\,.
\een
It is called the Mar\v cenko-Pastur distribution with ratio index $c$ and scale index $\sigma^2$ (\cite{Baimethodo} p.621).

It is well known (\cite{MarchPastur1}, \cite{Baimethodo} section 2.1.2) that as $n,r \rightarrow \infty$ with 
$r/n \rightarrow c \in (0, \infty)$, the family of empirical spectral distributions 
$(\mu_{n,r})$ 
 converges a.s. weakly
to the Mar\v cenko-Pastur distribution $\pi_1^c$. If we replaced the common law ${\cal N}(0 ; 1)$ by 
${\cal N}(0, \sigma^2)$ then the limiting distribution would be $\pi_{\sigma^2}^c$.

Recently, De Cock et al. (\cite{Spince1}) and Jiang (\cite{Jiang2}) proved that the same result holds true in the uniform Gram ensemble. 
If $\widetilde\lambda_{n,r}^{(k)}, k=1, \cdots , r$ be the (real) eigenvalues of
${\cal Y}_{n,r}$ and
\be
\widetilde\mu_{n,r} = \frac{1}{r} \sum_{k=1}^r \delta_{\widetilde\lambda_{n,r}^{(k)}}
\ee
then, as $n\rightarrow \infty$ and $r/n \rightarrow c \in (0, \infty)$, the family $(\widetilde\mu_{n,r})$ converges a.s. to $\pi_1^c$.

In both cases (uniform Gram and Wishart) we examine in  Section \ref{Hr} and 
 \ref{Wm} the connection between the "decomposition" method and the spectral method, at the level of marginals.

\section{Determinants in the uniform Gram ensemble}
The proofs of the results of this Section are in Section \ref{sect7}. The subscript or superscript G (resp. W) for the limiting quantities refers to the uniform Gram ensemble (resp. the Wishart ensemble). 

Let us notice that in the paper \cite{gamb}, 
 a decomposition in product of beta variables for a completely different problem leads to similar results.

\label{Hr}
\subsection{Two first moments and almost sure convergence}
\label{2first}

\begin{prop}
\label{first2}
For the mean, we have
\ben
\label{enough}
\lim_n \  \sup_{p \leq n} \left|  
 \frac{1}{n}`E\!\  \Upsilon_{n , p} + {\cal J}\left(1 - \frac{p}{n}\right)\right| = 0\,,
\een
and actually,
\ben
\label{espt}
\forall t \in [0,1) \ , \   \ \left(`E \!\  \Upsilon_{n}(t) + n{\cal J}\Big(1 - \frac{\lfloor nt\rfloor}{n}\Big)\right) 
&\rightarrow&  \d_G(t) := t + \frac{1}{2}\log (1-t)
\\
\label{esp1}
   \ \left(`E \!\  \Upsilon_{n,n} +n + \frac{1}{2} \log n\right) &\rightarrow& \frac{-\gamma - \log 2 +3}{2}
\een
where $\gamma$ is the Euler constant.\\
For the variance, we have 
\ben
\label{vart}
\forall t \in [0,1) \ , \  \  \ \hbox{Var} \  \Upsilon_n (t) &\rightarrow& \sigma^2_G (t) := -2\log (1-t) -2t\\
\label{var1}
   \ \left(\hbox {Var} \  \Upsilon_{n,n} - 2 \log n \right) &\rightarrow&  \frac{8\gamma + \pi^2}{4}\,.
\een
\end{prop}

\begin{theo}
\label{theocvps}
Almost surely, 
\begin{eqnarray*}
\lim_n \ \sup_{t\in [0,1]} \left|\frac{ \Upsilon_n (t)}{n} + {\cal J}(1-t)\right| = 0\,.
\end{eqnarray*}
\end{theo}
\medskip

The formulae (\ref{esp1}) and (\ref{var1}) are due to Abbott and Mulders (lemmas 4.2 and 4.4 in \cite{abbott}), using a variant of the decomposition method.

If we want to use the spectral method (with fixed $t=c < 1$) we may start with  
$$\frac{1}{n}\Upsilon_n (c) = \frac{\lfloor nc\rfloor}{n}\int \log x \ d\widetilde\mu_{n, \lfloor nc\rfloor}(x)$$
use the weak convergence of 
$\widetilde\mu_{n, \lfloor nc\rfloor}$ towards $\pi_1^c$,
(see Section \ref{MP}). To conclude that 
\ben
\label{logH}
\lim_n \int \log x \ d\widetilde\mu_{n, \lfloor nc\rfloor} (x)
= \int \log x \!\ d\pi_1^c(x)\,,
\een an additional control is necessary, since $x \mapsto \log x$ is not bounded. 
In \cite{Jiang2}, Jiang proved recently that the largest and the smallest eigenvalue 
of ${\cal Y}_{n,r}$ converge a.s. , as $r/n \rightarrow c < 1$ to $b(c)< \infty$ and $a(c)> 0$ respectively (remember the
definitions of $a$ and $b$ in (\ref{aetb})). So, (\ref{logH}) is true. Moreover,
 it is known  
(\cite{jonsson1} p.31 and \cite{BaiSilver2} p. 596-597) that :
\ben
\label{often}
c \int \log x  \ d\pi_1^c(x) &=& \int_{a(c)}^{b(c)}   \frac{\log x}{2\pi x}\ \sqrt{(x-a(c)) (b(c)-x)}\ dx\\
\nonumber&=& 
(c-1) \log(1-c) -c = -{\cal J}(1-c)
\een
This matches with the result of Theorem \ref{theocvps}.
\subsection{Fluctuations}
\label{flucH}
 Let $D_T = \{ v \in \BBd ([0,T]) : v(0) = 0\}$ and  $D= \{ v \in \BBd([0,1))  : v(0) = 0\}$
 the set of c\`adl\`ag functions on $[0,T]$ and $[0, 1)$, respectively, starting from $0$.
 
\begin{theo}
\label{DonskerH}
\begin{enumerate}
\item Let for $n \geq 1$
\begin{eqnarray*}
\eta_n^G (t) :=  \Upsilon_n (t)
 + n{\cal J}\Big(1 - \frac{\lfloor nt\rfloor}{n}\Big)
 \ \ , \ \ t\in [0,1)\,.
 \end{eqnarray*}
Then as $n \rightarrow \infty$
\ben
\label{etantH}
\left(\eta_n^G (t); \ t \in [0,1)\right) \Rightarrow \left( Y_t ; \ t \in [0,1)\right)\,,
\een
where $Y$ is the (Gaussian) diffusion solution of the stochastic differential equation :
\ben
\label{sdeG}
dY_t = \frac{1-2t}{2(1-t)}\!\ dt + \sqrt{\frac{2t}{1-t}}\ d{\bf B}_t\,,
\een
with $Y_0 = 0$, ${\bf B}$ is a standard Brownian motion and $\Rightarrow$ 
stands for the weak convergence of distributions in 
$D$ provided with the Skorokhod topology.

\item Let
\begin{eqnarray*}
\widehat\eta_n^G (1) = {\frac{\Upsilon_{n,n} +n + \frac{1}{2} \log n}{\sqrt{2\log n}}}\ \,.
\end{eqnarray*}
Then as $n \rightarrow \infty$, $\widehat\eta_n^G (1)\Rightarrow N$
where $N$ is a standard normal variable independent of $\bf B$, (and $\Rightarrow$ 
stands for the weak convergence of distribution in $\BBr$).
\end{enumerate}
\end{theo}

\subsection{Large deviations}
\label{ldH}
All along this section, as in Section \ref{ldpW} and in the proof Sections \ref{7.4} and \ref{similarroute}, 
we use the notation of Dembo-Zeitouni \cite{DZ}. In particular we write LDP for Large Deviation Principle. The reader may have some interest  in consulting
\cite{gamb}
where a similar method is used for a different model, but here
we use  a slightly different topology to be able to catch the marginals  in $T$.

For $T < 1$, let ${\cal M}_T$ be the set of  signed measures on $[0,T]$ and let ${\cal M}_<$ be the set of 
measures whose support is a compact subset of
 $[0,1)$. 

We provide $D$ with the weakened topology $\sigma(D, {\cal M}_<)$. 
So, $D$ is the projective limit of the family, indexed by $T< 1$ 
of topological spaces  $\big(D_T , \sigma(D_T,{\cal M}_T)\big)$.

Let $V_\ell$ (resp. $V_r$) be the space of left (resp. right)  continuous $\BBr$-valued functions with 
bounded variations. We put a superscript $T$ to
specify the functions on $[0,T]$. There is a one-to-one correspondence between $V_r^T$ and ${\cal M}_T$ : 

- for any $v \in V_r^T$, there exists a unique 
$\mu \in {\cal M}_T$ such that $v = \mu([0, \cdot])$; we denote it by $\dot v$

- for any $\mu \in {\cal M}_T$,  $v = \mu([0, \cdot])$ stands in $V_r$.
 
For the following statement, we need some notation. 
Let ${\bf H}$ be the entropy function : 
\ben{\bf H}(x|p) =\displaystyle x \log \frac{x}{p} + (1-x) \log\frac{1-x}{1-p}\,,\een
and set\footnote{we set $\delta(y |A) = 0$ if $y\in A$ and $= \infty$ if $y \notin A$}
\ben
\label{LaLs}
L_a^G (t, y) = \frac{1}{2}{\bf H}(1-t | e^{y}) \ \delta(y | (-\infty, 0))\  \ 
\ , \ L_s^G (t , y) = -\frac{1}{2}(1-t)y \ \delta(y | (-\infty, 0))
\,.
\een
\begin{theo}
\label{LDPH}
The sequence $\{
\frac{1}{n}\!\ 
 \Upsilon_n  \}_n$ satisfies a LDP in $(D, \sigma(D, {\cal M}_<))$ at scale $n^2$ 
with good rate function given for $v \in D$ by:
\ben
I_{[0,1)}^G (v) = \int_{[0,1)} L_a^G\left(t, \frac{d\dot v_a}{dt} (t)\right) dt + \int_{[0,1)} L_s^G \left(t, \frac{d\dot v_s}{d\mu} (t)\right) d\mu(t)\ \ \ if \ \ v \in V_r\,,
\een
where $\dot v = \dot v_a + \dot v_s$ 
is the Lebesgue decomposition of the measure $\dot v \in {\cal M}([0,1))$ in absolutely continuous and singular parts with respect 
to the Lebesgue measure and $\mu$ 
is any bounded positive measure dominating $\dot v_s$. If $v \notin V_r$, then $I_{[0,1)}^G (v) = \infty$.
\end{theo}

That means, roughly speaking, that 
$$`P( \Upsilon_n \simeq v) \approx \e^{-n^2 I_{[0,1)}^G (v)}\,.$$
The proof, in Section \ref{7.4} needs several steps. 
First we show that $\{\frac{1}{n}\!\ \dot \Upsilon_n\}$  satisfies a LDP in ${\cal M}_T$ provided with the topology 
$\sigma({\cal M}_T, V_\ell)$. Then we carry the LDP to  $\big(D_T , \sigma(D_T, {\cal M}_T)\big)$ with good rate function:
\ben
\label{42}
I_{[0,T]}^G (v) = \int_{[0,T]} L_a^G\left(t, \frac{d\dot v_a}{dt} (t)\right) dt + \int_{[0,T]} L_s^G 
\left(t, \frac{d\dot v_s}{d\mu} (t)\right) d\mu(t)\,.
\een
 To end the proof we apply the  
Dawson-G\"artner theorem on projective limits (\cite{DZ} Theorem 4.6.1, see also \cite{leo1} Proposition A2). 

Let us notice that $I_{[0,T]}^G$ vanishes only for $\frac{d\dot v_a}{dt}(t) = \log (1-t)$ and $\frac{d\dot v_s}{d\mu}(t) = 0$ (essentially) 
i.e. for $v (t) = - {\cal J}(1-t)$, which is consistent with the result of Theorem \ref{theocvps}.
\medskip

The LDP for marginals is given in the following theorem, where a rate function with affine part appears.
\begin{theo}
\label{margH}
For every $T < 1$, the sequence $\big\{\frac{1}{n}\!\  \Upsilon_{n, \lfloor nT\rfloor}\big\}_n$  satisfies a LDP in $\BBr$ 
at scale  $n^2$ with good rate function
denoted by 
\ben
\label{opt1}
I_T^{G} (\xi) = \inf \{I^G_{[0,T]}(v) \ ; \ v(T) = \xi\}\,.
\een
\begin{enumerate}
\item
 If $\xi \geq -T$ the equation
\ben
\label{=xi1H}
{\cal J}(1+2\theta) -{\cal J}(1-T+2\theta)  - T \log (1+2\theta) =\xi\,,
\een
has a unique solution, and we have
\ben
\label{idetH}
I_T^{G}(\xi) = \theta\xi &+& \frac{T}{2} {\cal J} (1+2\theta)\\
\nonumber
&+& \frac{1}{2}\left(F(1)-F(1-T) -F(1+2\theta) +F(1-T+2\theta)\right)
\,.
\een
\item
 If $\xi < - T$, we have 
\ben
\label{mauvH}
I_T^{G}(\xi) = I_T^{G}(-T) - \frac{1-T}{2}(\xi + T)\,.
\een
\end{enumerate}
\end{theo}

\section{Determinants in the Wishart ensemble}
\label{Wm}
\subsection{Introduction}
Three ways are possible to study the asymptotic behavior of the determinant of
$${\cal X}_{n,r} = \frac{1}{n}B'_{n,r}B_{n,r}
\,.$$

a) The spectral approach if we are intereted only in marginals ($r/n \rightarrow c$ fixed).

b) The Bartlett's decomposition method for a dynamical study. The representation (\ref{secondarray})  leads to results similar to those 
of the above section. Let us remark that at the level of marginals,
  (\ref{bart}) gives the 
Mellin transform:
\be
`E \left(\det {\cal X}_{r,n}\right)^s = 2^{rs} \frac{\Gamma_r \left(\frac{n}{2}+ s\right)}
{\Gamma_r \left(\frac{n}{2}\right)}\ , \ 
\hbox{for}\ \Re s > -\frac{n-r+1}{2}
\ee
where
$$\Gamma_r (\alpha) = \pi^{r(r-1)/4}\Gamma(\alpha)\Gamma\Big(\alpha-\frac{1}{2}\Big)\cdots \Gamma\Big(\alpha-\frac{r-1}{2}\Big)$$
(see for instance \cite{Mathai} Theorem 1 p. 347).
This yields
the density of $\det {\cal X}_{r,n}$ (cf. \cite{Mathaibook}, \cite{Mehta1} formula 2.4, when $r=n$), in which the 
Meijer function is involved. 

c) Actually, we prefer to establish these results from the representation (\ref{HW}) which we recall here:
\be
\log \det {\cal X}_{n,r} &=& \log \det {\cal Y}_{n,r} + S_{n,r} \\
\nonumber
S_{n,r}&=&\sum_{k=1}^r \log \frac{\Vert b_k\Vert^2}{n}\,,
\ee
for $r=1 , \cdots, n$, where the variables $\Vert b_k\Vert^2, k = 1, \cdots , n$ are independent,  $\chi^2_n$ distributed, 
and independent of $ (\log \det {\cal Y}_{n,r} =  \Upsilon_{n,r}, r= 1, \cdots, n)$.

We state also connections with known results deduced from the spectral approach.
The proofs are in Section \ref{sect7}.

\subsection{Almost sure convergence and fluctuations}
\label{ascfW}

By extension of the study in Section \ref{2first} (method b) above), we get the following.
\begin{prop}
\label{ODD}
For the mean we have
\be
\lim_n \sup_{p\leq n}\left|\frac{1}{n}\!\ `E\, \Xi_{n,p}
 + {\cal J}(1-\frac{p}{n})\right| = 0\,.
\ee
Moreover, as $n \rightarrow \infty$
\be
\forall t \in [0,1) \ \ \ `E
\, \Xi_{n, \lfloor nt\rfloor}
+ n{\cal J}\Big(1 - \frac{\lfloor nt\rfloor}{n}\Big) 
&\rightarrow&   \d_W(t) = \frac{1}{2}\log (1-t)
\\
 `E\,   
 \Xi_{n,n}
  +n + \frac{1}{2} \log n &\rightarrow& -\frac{\gamma + \log 2 -1}{2}\,.
\ee
For the variance we have
\be
\forall t \in [0,1) \ \ \ \hbox{Var}\,  
\Xi_{n, \lfloor nt\rfloor}
&\rightarrow& \sigma^2_W(t) := - 2\log (1-t)\\
\hbox{Var}\,  \Xi_{n,n}
 - 2 \log n   &\rightarrow& \frac{2\gamma + \pi^2 +8}{4}\,.
\ee
\end{prop}

\begin{theo}
\label{theocvpsW}
Almost surely, 
\be
\lim_n \ \sup_{t\in [0,1]} \left|
\frac{1}{n}\!\ \Xi_{n, \lfloor nt\rfloor}  + {\cal J}(1-t)\right| = 0\,.
\ee
\end{theo}

\begin{theo}
\label{Donsker}
Let
\begin{eqnarray*}
\eta_n^W (t) &:=& 
\Xi_{\lfloor nt\rfloor,n} + n {\cal J}\Big(1 - \frac{\lfloor nt\rfloor}{n}\Big) 
 \ \ , \ \ t\in [0,1]\,,\\ 
\widehat\eta_n^W (1) &=& \frac{
\Xi_{n,n} +n + \frac{1}{2}\log n}{\sqrt{2\log n}}\ \,.
\end{eqnarray*}
Then as $n \rightarrow \infty$
\ben
\label{etant}
\left(\eta_n^W (t); \ t \in [0,1)\right) &\Rightarrow& (X_t, \ \ t \in [0,1))\\
\nonumber
\widehat\eta_n^W (1)&\Rightarrow& N
\een
where $X$ is the Gaussian diffusion solution of the stochastic differential equation :
\ben
\label{sdeW}
d X_t = - \frac{1}{2(1-t)}\!\ dt + \sqrt{\frac{2}{1-t}}\ d {\bf B}_t\,,
\een
with $X_0=0$, where
${\bf B}$ is a standard Brownian motion and $N$ is a standard normal variable independent of ${\bf B}$.
\end{theo}

With the direct method b), Jonsson  proved (\ref{etant}) for fixed $t$ 
(i.e. convergence in distribution of the marginal) and deduced a convergence in probability 
of $\frac{1}{n}
\Xi_{n, \lfloor nt\rfloor}$  towards $-{\cal J}(1-t)$ (Theorem 5.1 p.29 and Corollary 5.1 p.30 of \cite{jonsson1}).

 Let us explain now the results which may be obtained by the spectral method a). 

It is well known that the empirical spectral distribution of ${\cal X}_{n,r}$ 
converges weakly a.s. when $r/n \rightarrow c$ towards $\pi_{1}^c$.  Moreover if $c < 1$ the largest (resp. smallest) 
eigenvalue converges a.s. to $b(c)$ (resp. to $a(c)$). For comments on these results and references, one may consult \cite{Baimethodo} 
sections 2.1.2 and 2.2.2. 

In our context, this implies easily that 
\ben
\label{esp}
\frac{1}{n}\, \Xi_{n,r} = \frac{1}{n}\log \det {\cal X}_{n, r} = \frac{r}{n}\int \log x \ d\mu_{n,r}(x) \rightarrow c \int \log x  \ d\pi_1^c(x)
\een
a.s. when $r/n \rightarrow c \in (0,1)$. We already see in (\ref{often}) the value of this integral. 
Claim (\ref{esp}) is then consistent with Theorem \ref{theocvps}.

The fluctuations were studied recently by Bai and Silverstein \cite{BaiSilver2} (in the  case of complex entries). They obtained 
\be
\log \det {\cal X}_{n,r}  
+ 
n{\cal J}\Big(1 - \frac{r}{n}\Big)
\Rightarrow {\cal N}\left(\log (1-c)\ 
; \ -2\log (1-c)\right)\,, 
\ee
which is consistent with (the marginal version of) (\ref{etant}). 
\medskip
\subsection{Large deviations}
\label{ldpW}
Again, the three routes are possible to tackle the problem of large deviations for determinant of Wishart matrices. 
A direct method would use the cumulant generating function from Section \ref{bartlettr} and would meet computations 
similar to those seen in the Gram case.

 To avoid repetitions, we use the b) method, drawing benefit from  
an auxiliary study of  $S_{n,r}$.
\begin{lem}
\label{LDPaux}
The 
sequence  $\big\{ \frac{1}{n}\!\ S_{n}(t) , 
t \in [0, 1)\big\}_n$ satisfies a LDP in $(D, \sigma(D, {\cal M}_<))$
at scale $n^2$ with good rate function
\be
\label{59}
I_{[0,1)}^{S}(v) = 
\int_{[0,1)} L_a^{S}\left(\frac{d\dot v_a}{dt} (t)\right) dt + \int_{[0,1)} L_s^{S} \left(\frac{d\dot v_s}{d\mu} (t)
\right) d\mu(t)
\ee
where 
\ben
\label{60}
L_a^{S} (y) = \frac{1}{2}\left(e^{y} -y-1\right) \ , \ L_s^S (y) = -\frac{y}{2}\delta(y|(-\infty,0))\,,
\een
and  $\mu$ is any  measure dominating $d\dot v_s$. 
\end{lem}
Let us stress that the instantaneous rate functions  are time homogeneous 
and then we may write $[0,1]$ instead of $[0,1)$.
\begin{theo}
\label{LDPwish}
The 
sequence $\big\{\frac{1}{n}\!\ 
\Xi_{n, \lfloor nt\rfloor}, t \in [0, 1)\big\}_n$
satisfies a LDP in $(D, \sigma(D, {\cal M}_<))$, at scale $n^2$ with good rate function 
\ben
\label{a+s}
I^{W}_{[0,1)} (v) = \int_{[0,1)} L_a^{W}\left(t, \frac{d\dot v_a}{dt} (t)\right) dt + 
\int_{[0,1)} L^{W}_s \left(t, \frac{d\dot v_s}{d\mu} (t)\right) d\mu(t)
\een
where
\ben
\label{identif}
\nonumber
L_a^{W}(t,y) &=&
 \frac{1}{2}(e^y -1)  - \frac{1}{2} (1-t)y +\frac{1}{2}{\cal J}(1-t)\\
L_s^{W} (t,y) &=&  -\frac{1}{2} (1-t) y \ \delta(y | (-\infty , 0))
\,.
\een 
and  $\mu$ is any  measure dominating $d\dot v_s$. 
\end{theo}
Let us notice that the restriction of $I_{[0,1)}^W$ to $[0,T]$ vanishes only for $\frac{d\dot v_a}{dt}(t) = \log (1-t)$ and $\frac{d\dot v_s}{d\mu}(t) = 0$ (essentially) 
i.e. for $v (t) = - {\cal J}(1-t)$, which agrees with the result of Theorem \ref{theocvpsW}.
\medskip

For marginals,  we give (without proof) the exact analogue of Theorem \ref{margH}.
\begin{theo}
For every $T < 1$, the sequence $\big\{\frac{1}{n}\!\ 
\Xi_{n, \lfloor nT\rfloor}\big\}_n$  satisfies a LDP 
in $\BBr$ at scale  $n^2$ with good rate function
denoted by $I_T^{W}$.
\begin{enumerate}
\item If $\xi \geq \xi_T := {\cal J}(T) - 1$ the equation
\ben
\label{=xi1W}
 {\cal J}(1+2\theta) -{\cal J}(1-T+2\theta) =\xi\,.
\een
has a unique solution, and we have
\ben
\label{idetW}
I_T^{W}(\xi) =  \theta\xi 
+ \frac{1}{2}\left(F(1)-F(1-T) -F(1+2\theta) +F(1-T+2\theta)\right) 
\,.
\een
\item
 If $\xi < \xi_T$
, we have 
\be
I_T^{W}(\xi) = I_T^{W}(T) + \frac{1-T}{2}(\xi_T-\xi)\,.
\ee
\end{enumerate}
\end{theo}
\medskip

Let us comment the spectral approach.
Hiai et Petz \cite{hiai1} proved that if $n \rightarrow \infty, r/n \rightarrow c < 1$, then $\{\mu_{n,r}\}$ satisfies 
a LDP at scale $n^2$ with some explicit good rate function $I^{(sp)}_c$ given below in (\ref{HP1}, \ref{HP2}, \ref{HP3}). 
If the contraction $\mu \mapsto \int \log x \ d\mu(x)$ 
 was continuous, we would claim that
$\{\frac{1}{n} \log \det {\cal X}_{n,\lfloor nT\rfloor}\}_n$ satisfies a LDP with good rate function
\be
I^{W}_T (\xi) = \inf \left\{I^{(sp)}_T (\mu) \ ; \ T\int \log x 
\ d\mu(x) = \xi \right\}\,.
\ee 
Actually
\ben
\label{HP1}
I^{(sp)}_c (\mu) := -\frac{c^2
}{2} \Sigma(\mu) + \frac{c
}{2} \int\Big(x - (1-c) \log x\Big) d\mu(x) + B(c)
\een
where
\ben
\label{HP2}
\Sigma(\mu) := \int\!\!\int \log|x-y|\ d\mu(x) d\mu(y)
\een
is the so-called logarithmic entropy and for $c \in (0,1)$
\ben
\label{HP3}
B(c) = -\frac{1}{4}\Big( 3c -c^2\log c + (1-c)^2 \log(1-c)\Big)\,.
\een
We do not know if the contraction 
$\mu \mapsto \int \log x \ d\mu(x)$ 
works, although not continuous. However we will prove the following result, where for $u \in `R$ we set ${\cal A}(u) = \{\mu : \int \log x \ d\mu(x) =u\}$. 
\begin{prop}
\label{infI} 
For $\xi\geq \xi_T$ and $\theta$ solution of (\ref{=xi1W}), let  $\sigma^2 = 1 +2\theta$.
Then the infimum of $I_T^{(sp)} (\mu)$ on ${\cal A}(\xi/T)$ is uniquely achieved for $\pi_{\sigma^2}^{T/\sigma^2}$
 and 
\ben
\label{defxi} 
I_T^{W}(\xi) = I_T^{(sp)} (\pi_{\sigma^2}^{T/\sigma^2})
 = \inf \{ I_T^{(sp)} (\mu) ; \ \mu \in {\cal A}(\xi/T)\}\,.
\een
\end{prop}

\bigskip

\begin{rem} 
\begin{enumerate}
\item The endpoint is $\xi_T = {\cal J}(T) - 1$, corresponding
 to $1+2\theta = T$, 
i.e. $\sigma^2 = T$.

\item For $\xi < \xi_T$ we do not know what happens. 
We can imagine that the infimum  in (\ref{defxi}) has a solution in some extended space.
\end{enumerate}
\end{rem}
\section{Extensions}\label{Je}
We examine now some possible extensions of the previous results. We focus on  assumptions on the entries of the matrix $B$. 
We keep the same asymptotics and notation than in above sections.

Let us mention that the methodology of the present paper will be applied to the Jacobi ensemble in a forthcoming paper.

\subsection{Independent Gaussian non real entries}
In the previous sections, entries of $B$ were real numbers. We consider now entries in $\BBc$ and in $\BBh$, 
the set of real quaternions.
Recall that an element of $\BBh$ may be viewed as 
a $2\times 2$ matrix of the form
\begin{eqnarray*} 
q= \left( \begin{array}{cc}
z & w \\
-\bar w& \bar z \end{array}\right)
\end{eqnarray*}
where $z$ and $w$ are complex numbers. 	  Its dual (or conjugate) is
\begin{eqnarray*} 
\bar q = \left( \begin{array}{cc}
\bar z & -w \\
\bar w&  z \end{array}\right)\,.
\end{eqnarray*}
denoted by $z_{jk}$  for $j = 1, \cdots ,n$ and $k=1, \cdots , r$.

Let $B$ be a $n\times r$ random matrix, and suppose the entries of $B$ are determined by a parameter $\beta = 1, 2$ or $4$. These entries are i.i.d. random variables $\BBr, \BBc$ or $\BBh$ valued 
 with Gaussian densities
\be
\frac{1}{\sqrt{2\pi}}\!\ \e^{- b_{jk}^2 /2}\ , \ \ \frac{1}{{\pi}}\!\ \e^{- |z_{jk}^2|}\ , 
\ \ \frac{2}{{\pi}}\!\ \e^{- 2|z_{jk}^2| }\ , \ \ \frac{2}{{\pi}}\!\ \e^{- 2|w_{jk}^2|}\,.
\ee 
in the three cases $\beta = 1,2$ and $4$ respectively.
If 
$B^\dag$ denotes the transpose of the conjugate of $B$, then the $r\times r$ matrix ${\cal X}
 = B^\dag B$  belongs to the Laguerre orthogonal (resp. unitary, resp. symplectic) ensemble  denoted 
 LOE (resp. LUE, resp. LSE).
 
Two main features of the LOE are shared by the LUE and LSE.  

$\bullet$  The Barlett decomposition still holds. For $\beta =2$ 
the references are for instance \cite{Goodman}, \cite{Letac03} (see also \cite{Forrester} Proposition 2.12 and 2.13). By the same argument as in section 
\ref{bartlettr}, 
the random variables  $(R_{k,j} , \ k=1 , \cdots, j-1)$ 
are complex i.i.d. normal of variance $1$, hence $2|R_{k,j}|^2  \buildrel{\cal D}\over{=} \chi_{2}^2$. Since 
$2\Vert b_j\Vert^2  \buildrel{\cal D}\over{=} \chi_{2n}^2$
we conclude that $2 R_{j,j}^2  \buildrel{\cal D}\over{=} \chi_{2n-2(j-1)}^2$.
For $\beta = 4$, the reference is \cite{Forrester} Exercise 2.5.5. We conclude that 
$4 R_{j,j}^2  \buildrel{\cal D}\over{=} \chi_{4n-4(j-1)}^2$.

A pathwise study of $\log \det {\cal X}$ in the cases $\beta =2$ or $4$ needs only 
slight modifications of arguments and  would lead to results very similar to those of Section \ref{Wm}.

$\bullet$ The spectral approach is built on the probability density of the eigenvalues 
$\lambda_j , j = 1, \cdots , r$ of ${\cal X}$ which is proportional to
\be
\prod_{j=1}^r \lambda_j^{\frac{\beta(n-r + 1)}{2} -1} \e^{-\beta\lambda_j/2}\ \prod_{1 \leq j < k \leq m} |\lambda_k - \lambda_j|^\beta
\ee
Convergence of the ESD is known not only for $\beta = 2, 4$ but for every $\beta > 0$ 
(see for instance the Dumitriu thesis \cite{Dum1} Theorem 6.5.1). 

Besides, the large deviations treated in \cite{hiai1} are stated for the (real) Wishart ensemble, 
but of course are available in the general case with slight modifications since everything rests on their Theorem 1. 

\subsection{Independent isotropic columns} We keep independence of vectors $b_1, \cdots , b_n$ but assume only isotropy (in $`R^n$) 
of their common distribution $\nu_n$. The polar decomposition allows  to obtain similar results as in Section \ref{Wm} under convenient assumptions on the radial distribution. 
 Let $\varepsilon_n = \log \Vert b_1 \Vert^2 - \log `E \Vert b_1 \Vert^2$ (remember that we omit the dimension index $n$). 
 
 To get convergence and fluctuations it is enough to assume 
 \ben
 \label{hypr}
 n `E \varepsilon_n \rightarrow a_1 \ , \ n\hbox{Var}\, \varepsilon_n \rightarrow a_2 \ , \ n `E (\varepsilon_n - `E \varepsilon_n)^4 \rightarrow 0\,.
 \een
To get large deviations, it would be sufficient to assume that 
$\frac{1}{n^2} \sum_{k=1}^n \ \log `E \exp \big(n\varphi(k/n)\varepsilon_n\big) $  
has a limit for some convenient functions $\varphi$.

Notice that in \cite{akhavi2} and \cite{akhavi3}, the authors 
use the uniform distribution 
in the unit {\it ball}, so that the distribution of $\Vert b_1\Vert^2$ is beta$\big(\frac{n}{2}, 1\big)$ and (\ref{hypr}) 
is satisfied with 
$a_1 = -2 , \ a_2 = 0$. Here the contribution of the radial part is roughly "deterministic" since $`E \Vert b_1\Vert^2$ is bounded.

\subsection{Independent identically distributed (non Gaussian) entries}
If we restrict ourselves to marginals only, we may leave the Gaussian world.  Let us 
assume i.i.d. (real) entries with finite variance. 
In 
(\cite{jonsson1}, \cite{BaiSilver}), the authors 
 proved the 
convergence of the spectral distribution,
using Stielj\`es' transform (\cite{MarchPastur1}). 
In \cite{BaiSilver2} Bai and Silverstein assumed 
 $`E b_{11}^4 = 3$ (real entries) or 
$`E b_{11}^2 = 0$ and $`E b_{11}^4 = 2$ (complex entries), and proved a central limit theorem for linear statistics, with 
the meaningful example of logarithm of the determinant.

The Bartlett's decomposition is not possible in the general case. Nevertheless, a product formula for the determinant is well known (see for example  Lemma 3.1 p.9  and formula 4.3 p.15 in \cite{Zeit}) but nothing can be said about 
the distribution of the components of the product.

Moreover, using again the norming of column vectors as in previous sections, we may define $\widetilde B$ with
$$\widetilde b_{ij} = \frac{b_{ij}}{\sqrt{\sum_{k=1}^n b_{kj}^2}}$$
 A slight modification of this matrix is used in multivariate analysis to test  that variates are uncorrelated. The matrix 
 ${\cal Y} = \widetilde B' \widetilde B$ is a Gram matrix, built from  
 independent vectors, identically distributed and living on $`S_{n-1}$. 
 In this context, Jiang (\cite{Jiang2}) recently proved the convergence of the ESD to $\pi_1^c$ and also the convergence of the extreme eigenvalues. It is then easy to 
deduce the convergence of the normalized logarithm of the determinant. 

This model is clearly an extension of the uniform Gram ensemble, for which De Conck et al. (\cite{Spince2}) 
proved the convergence of the ESD to $\pi_1^c$ 
 with an independent method. 

\subsection{Independent isotropic rows}  We keep independence of rows of the random matrix $B$ and 
assume that they are identically distributed with an isotropic distribution 
on $`R^r$. Actually in the data matrix $B_{n,r}$, the index $n$ 
is, as previously, the size of the sample and $r$ is the number of variates.  In \cite{Yin1}, 
Yin and Krishnaiah  proved the convergence of the ESD of ${\cal X}_{n,r} = \frac{1}{n} B_{n,r}'B_{n,r}$ but 
 the limiting distribution was not known. Actually, when the underlying distribution is uniform on $`S_{r-1}$ 
  we can identify the  limiting distribution from the result of Jiang \cite{Jiang2} or De Cock et al. \cite{Spince1}.

 We set $C = B'$ and then  $C'C$
 is in the uniform Gram ensemble (in $M_{nn}$). The eigenvalues of $C'C$ 
are (except $0$ with multiplicity $n-r$) the same as those of $CC' = B'B$. If $\mu^\star_{n,r}$ is the ESD of $CC'$, 
the ESD of $C'C$ 
is then $$\mu_{n,r}= \frac{r}{n}\mu^\star_{n,r} + (1- \frac{r}{n})\delta_0$$
If $r/n \rightarrow c < 1$, hence $n/r \rightarrow 1/c > 1$,
$$\lim_{n,r \rightarrow \infty}\mu_{n,r}  = \pi_1^{1/c}$$
(\cite{Spince1} Theorem 10 or \cite{Jiang2} Theorem 2)
so that, $$\lim_{n,r \rightarrow \infty}\mu^\star_{n,r}  = \frac{1}{c}\left(\pi_1^{1/c} - (1-c)\delta_0\right)$$ and 
 from (\ref{defpi}) we see that
\be
\lim_{n,r \rightarrow \infty}\mu^\star_{n,r}  
= \pi_{1/c}^c
\ee
Besides, Yin and Krisnaiah scaling is $\frac{1}{n}$, 
so that the limit of the ESD of ${\cal X}_{n,r}$ is the image of $\pi_{1/c}^c$ by the dilatation 
$D_c : x \mapsto cx$ i.e. $\pi_1^c$ (see Section \ref{MP}). Moreover, the results on extreme eigenvalues obtained by Jiang are easily carried.

\section{Proofs of Section \ref{Hr}}
\label{sect7}
To stress the dependence on $n$ we set for $2 \leq k\leq n$
$$h_{k,n} = \widetilde R_{k,k}^2\,.$$
(and $h_{1,n} = 1$). The key tool is the cumulant generating function :
\ben
\label{cgfH}
\Lambda_{n,k}^G (\theta) := \log \BBe \exp\left(\theta\log h_{k,n}\right) = \log \BBe \big(h_{k,n}^{\theta}\big)\,.
\een
From Proposition \ref{celebrated} 2), we have for $2 \leq k\leq n$ 
\ben
\label{LambdaHH}
\Lambda_{n,k}^G (\theta) = \ell\left(\frac{n-k+1}{2}+ \theta\right) - \ell\left(\frac{n-k+1}{2}\right)
- \ell\left(\frac{n}{2} +\theta\right) + \ell\left(\frac{n}{2}\right)
\een
where $\ell = \log \Gamma$.
By derivation
\be
`E\left( \log h_{k,n} \right) &=& \Psi\left(\frac{n-k+1}{2}\right) - \Psi\left(\frac{n}{2}\right) \\
\nonumber
\hbox{Var}\! \left( \log h_{k,n} \right) &=& \Psi'\left(\frac{n-k+1}{2}\right) - \Psi'\left(\frac{n}{2}\right) \,,
\ee
where $\Psi = \ell' = \Gamma'/\Gamma$ is the digamma function. 

\subsection{Proof of Proposition \ref{first2}}
We need the following lemma. 
\begin{lem}
\label{premier}
For every $p\leq n$, we have
\ben
\nonumber
\BBe\!\  \Upsilon_{n, p}
&=& \frac{n-1}{2}\Psi\left(\frac{n+1}{2}\right) + \left(\frac{n-2p}{2}\right)\Psi\left(\frac{n}{2}\right) + 1-p\\ 
\label{edH}
&-&\frac{n-p-1}{2}\Psi\left(\frac{n-p+1}{2}\right) -\frac{n-p}{2}\Psi\left(\frac{n-p+2}{2}\right) 
\een
and
\ben
\label{svarH}
\left| \hbox{Var}\,  \Upsilon_{n,p} - 2(H_n -H_{n-p} -\frac{p}{n} )
\right| 
 \leq 4 \sum_{k=n-p+1}^n \frac{1}{k^2}\,.\een
\end{lem}
\medskip

\noindent{\bf Proof of Lemma \ref{premier}:} 
From (\ref{defXi}) and (\ref{LambdaHH}) we get
\ben
\nonumber
\BBe \Upsilon_{n, p} = \sum_2^{p} \left(\Lambda_{n,k}^G\right)' (0) 
&=& \sum_2^{p} \left[\Psi\left(\frac{n-k+1}{2}\right) - \Psi\left(\frac{n}{2}\right) \right]\\
\label{lespsi}
&=& \sum_{k=n-p+1}^n \Psi\left(\frac{k}{2}\right) - p \Psi\left(\frac{n}{2}\right)
\,.
\een
From the classical identity
$$\Psi\left(\frac{k+2}{2}\right) - \Psi\left(\frac{k}{2}\right) = \frac{2}{k}$$
and since $\Psi(1/2) = -\gamma -2\log 2$, Abbott and Mulders \cite{abbott} deduced
\ben
\label{abbmuld}
\sum_{i=1}^{k-1}\Psi\left(\frac{i}{2}\right) = \frac{k-2}{2}\Psi\left(\frac{k}{2}\right) -k+ 
\frac{k-1}{2}\Psi\left(\frac{k+1}{2}\right)
+ \frac{2 -\gamma -2\log 2}{2}\,.
\een 
It remains to take successively $k = n+1$ and $k=n-p+1$ and use (\ref{lespsi}).  
Besides
$$\hbox{Var}\, \Upsilon_{n,p} 
= \sum_{k=2}^{p} \left(\Lambda^G_{k,n}\right)'' (0)= \sum_{k=2}^{p} 
\left[\Psi'\left(\frac{n-k+1}{2}\right) - \Psi'\left(\frac{n}{2}\right)\right]\,,$$
so that (\ref{svarH}) comes from (\ref{polygamma}) and from (\ref{restepsi}) with $q=2$.
\QED

\medskip
\noindent{\bf Proof of Proposition \ref{first2}:}
We have only to prove (\ref{enough}), (\ref{espt}) and (\ref{vart}).

1) We have  
$$n{\cal J}\left(1 - \frac{p}{n}\right) = p+ (n-p)\log (n-p) - (n-p)\log n\,.$$
for $p< n$ and $=n$ for $p=n$. 
Using (\ref{edH}) and (\ref{supx}) we get (\ref{enough}).

2) Now, for $p =\lfloor nt\rfloor$, $0<t<1$ and $n \rightarrow \infty$ we use the more precise estimate (\ref{supx2}) in (\ref{edH}).
We leave the details to the reader.

3) For the variance, we start from (\ref{svarH}) and we get easily   (\ref{vart}).
\QED

\subsection{ Proof of Theorem \ref{theocvps}}
Since ${\cal J}$ is uniformly continuous on $[0,1]$ we have
$$\lim_n \sup_{t\in [0,1]} \left|{\cal J}\left(1 - \frac{\lfloor nt\rfloor}{n}\right) - {\cal J}(1-t)\right| = 0\,,$$
so that, owing to (\ref{enough}), it is enough to prove that a.s.
$ \sup_{1\leq p\leq n}
 \left|\Upsilon_{p , n}  
- ` E\Upsilon_{p , n}\right| = o(n)$. 
Actually this convergence is a consequence of  Borel-Cantelli's lemma, Doob's inequality and of 
the variance estimate Var $\frac{1}{n}  \Upsilon_{n, n} = O(n^{-2}\log n)$ coming from 
(\ref{var1}).
\QED
\subsection{Proof of Theorem \ref{DonskerH}}
Let us first notice that, thanks to the estimations of expectations in (\ref{espt}) and (\ref{esp1}), 
we can reduce the problem 
to the centered process
$\Delta_n (t) := \Upsilon_n (t) - `E\!\ \Upsilon_n (t)$ 
and to the centered variable
 $\widehat\Delta_n = \Delta_n (1)/ \sqrt{\log n}$.
\medskip

\noindent 1) We use the notation of (\ref{espt}) and (\ref{vart}). 
We have $\Delta_n (t) = \sum_{k=1}^{\lfloor nt\rfloor} \eta_{n,k}$
where
\[\eta_{n,k} := ( \log h_{k,n}) - `E( \log h_{k,n}),\ \  k \leq n\] is a row-wise independent arrow. 
To prove (\ref{etantH}) it is enough to prove the convergence in distribution 
in $\BBd([0,T])$, for every $T < 1$, of $\Delta_n$ to a centered Gaussian process with independent increments, 
and variance $\sigma^2_G$.
To this purpose we  apply a version of the Lindeberg-L\'evy-Lyapunov criteria
(see \cite{DDC} Theorem 7.4.28 of the french edition, or \cite{JacShi} §3c). 
For $t < 1$, from (\ref{vart})
it is enough to prove that 
\ben
\label{ddc4H}
\lim_n \sum_{k=1}^{\lfloor nt\rfloor} `E\!\ (\eta_{n,k}^4)
= 0\,.
\een
 We have from (\ref{cgfH}) 
 \ben \label{4=1+3H} `E (\eta_{n,k}^4) 
 =
(\Lambda_{k,n}^G) ^{(4)} (0) + 3 [(\Lambda_{k,n}^G) '' (0)]^2\,. \een
On the one hand, from (\ref{LambdaHH}), (\ref{polygamma}) and (\ref{restepsi}) for $q=4$ 
\ben
\label{channel4H}
\Big|\sum_{k=1}^p (\Lambda_{k,n}^G) ^{(4)} (0) - 48
\sum_{k=1}^p\Big[ \frac{1}{(n-k + 1)^3} - \frac{1}{n^3}\Big] 
\Big)\Big|
 \leq 96 \sum_{k=1}^p \frac{1}{(n-k+1)^4}\,.
\een
which, for $0< t < 1$ and $p=\lfloor nt\rfloor$ yields 
$\lim_n \sum_{k=1}^{\lfloor nt\rfloor} 
(\Lambda_{k,n}^G) ^{(4)} (0) = 0\,.$ 
On the other hand, 
\ben \label{supH} \sum_{k=1}^p [(\Lambda_{k,n}^G) '' (0)]^2
\leq \big(\sup_{j\leq p}(\Lambda_{j,n}^G) '' (0)\big) \sum_{k=1}^p
(\Lambda_{k,n}^G) '' (0) \een 
and from (\ref{vart}) we get
 $\lim_n \sum_{k=1}^{\lfloor nt\rfloor} (\Lambda_{k,n}^G) '' (0) = \sigma^2_G(t)$.
Besides, $\Psi'$ is non-increasing (see (\ref{polygamma})) so that
\[\sup_{j\leq p}(\Lambda_{j,n}^G) '' (0)\leq  \Psi'\left(\frac{n-p+1}{2}\right) \]
and from (\ref{restepsi}) (again), this term  tends to $0$. We just checked
(\ref{ddc4H}), which proves that the sequence of processes $(\Delta_n(t), t \in [0, 1))_n$ 
converges to a Gaussian centered process  ${\cal W}$ with independent increments 
and variance $\sigma^2_G$. So, by (\ref{espt}),  $\eta_n^G$ converges to the Gaussian process 
${\cal W} +\d_G$ with independent increments, drift $\d_G$ and variance $\sigma^2_G$.

Finally, equation (\ref{sdeG}) comes from
$$\d_G(t) = \int_0^T \frac{1-2s}{2(1-s)} \!\ ds \ \ , \ \ \sigma^2_G(t) = \int_0^t \frac{2s}{1-s} \!\ ds\,.$$

\medskip

\noindent 2)  
When $t =1$, most of the sums studied above explode and we need a renormalisation.
In fact, for every $n$, the process $\big(\Delta_n (t), t\in [0,1]\big)$ has independent increments. 
The conditional distribution of 
$\Delta_n (1)$, 
knowing
  $\Delta_n (t_1)=\varepsilon_1, \cdots , \Delta_n (t_r) = \varepsilon_r$ 
for $t_1 < \cdots < t_r$ is the same as 
$\varepsilon_r + \sum_{[nt_r]+1}^n \eta_{k,n}$. 
Formulae (\ref{svarH}) and (\ref{var1}) yield
\be
\sum_{[nt_r]+1}^n `E (\eta_{k,n}^2) = 2 \log n + O(1)\,.
\ee
Actually we can apply the Lindeberg's theorem (with the criterion of Lyapunov)
to the triangular array of random variables $\widehat \eta_{k,n} = \eta_{k,n}/ {\sqrt{2\log n}}$ with 
with $k = [nt_r]+1, \cdots , n$.
It is enough to prove
\ben
\label{lyaH}
\lim _n \sum_{k=1}^n `E (\widehat \eta_{k,n}^4) = 0\,.
\een
The route is the same as before, starting from (\ref{4=1+3H}), but now 
(\ref{channel4H})
 says that the sum
  $\sum_{k=1}^n \Lambda_{k,n} ^{(4)} (0)$
 is bounded.
In (\ref{supH}), the sum on the right (with $p=n$) is now equivalent to $ 2\log n$ and
the supremum (with $p=n$) is bounded. This yields
\[
\sum_{k=1}^n `E (\widehat \eta_{k,n}^4) = (\log n)^{-2} \sum_{k=1}^n
`E (\eta_{k,n}^4) = O((\log n)^{-1} )
\]
which proves
 (\ref{lyaH}). 

Then $\sum_{[nt_r]+1}^n  \eta_{k,n}/ {\sqrt{2\log n}}$ 
converges in distribution to ${\cal N}(0,1)$, and the same is true for the conditional distribution of 
$\widehat\Delta_n$ knowing $\Delta_n (t_1)=\varepsilon_1, \cdots , \Delta_n (t_r) = \varepsilon_r$. 
Since the limiting distribution does not depend on $\varepsilon_1, \cdots , \varepsilon_r$, we have proved that 
$\widehat\Delta_n$ converges in distribution to a random variable which is 
${\cal N} (0,1)$ and independent of ${\cal W}$.
\QED
\subsection{Proof of Theorem \ref{LDPH}}
\label{7.4}
As mentioned after the statement of the theorem, we have to prove the LDP for the restriction of 
$\frac{1}{n}\!\ \dot\Upsilon_n$ to $[0,T]$, viewed as an element of ${\cal M}_T$, at scale $n^2$ 
with rate function
\ben
\label{tauxHT}
\widetilde I_{[0,T]}^G (m) := \int_0^T L_a^G\left(t, \frac{dm_a}{dt} (t)\right) dt + \int_0^T L_s^G 
\left(t, \frac{dm_s}{d\mu} (t)\right) d\mu(t)\,.
\een
Let $V_\ell$ be the set of functions from $[0,T]$ to $\BBr$ which are left continuous and have bounded variation, 
and let $V_\ell^{*}$ be its topological dual when $V_\ell$ is provided with the uniform convergence topology.
  
 Actually  $\frac{1}{n}\!\ \dot \Upsilon_n \in {\cal M}_T$ may be identified with an element of $V_\ell^{*}$ 
 (see \cite{leo1} Appendix B):
 its action on $\varphi \in V_\ell$ is given by
\be
<\frac{1}{n}\!\ \dot \Upsilon^n, \varphi> := \frac{1}{n}
\sum_{k=1}^{\lfloor nT\rfloor} \varphi\Big(\frac{k}{n}\Big)  \log h_{n,k}\,.
\ee

The proof is based on the ideas of Baldi's theorem (\cite{DZ} p.?). 
The main tool 
is the normalized cumulant generated function (n.c.g.f.) which here takes the form
\be
{\cal L}_{n, \lfloor nt\rfloor}^G (\varphi) := \frac{1}{n^2}\log `E \exp n<\dot \Upsilon^n , \varphi>
\ee
Owing to (\ref{cgfH}) we have
\ben
\label{cgfnG}
{\cal L}_{n, \lfloor nt\rfloor}^G (\varphi) = \frac{1}{n^2} \sum_{k=1}^{\lfloor nt\rfloor} 
\Lambda_{n,k}^G\left(n\varphi\left(\frac{k}{n}\right)\right)
\een
and from (\ref{LambdaHH}) it is finite if
$\varphi\left(\frac{k}{n}\right) > - \frac{n-k+1}{2n} $
for  $1 \leq k \leq \lfloor nT\rfloor$ and $+\infty$ otherwise.

In Subsection \ref{cvncgf}, we  prove the convergence of this sequence of n.c.g.f. for a large 
class of functions $\varphi$. It will be sufficient, jointly to the variational formula given 
in Subsection \ref{variaf} to get the upperbound for probability of compact sets. 
Then Subsection \ref{expotigh} is devoted to exponential tightness, which  allow to 
get the upperbound for closed sets. However, since the limiting n.c.g.f. is not defined everywhere, 
the lowerbound (for open sets) is more delicate than in Baldi's theorem. 
Actually a careful study of exposed points as in \cite{grz} is managed in Subsection \ref{exposed}. We end the proof in \ref{end}.

\subsubsection{Convergence of the n.c.g.f.}
\label{cvncgf}
\begin{lem}
\label{cvcgfH}
If $\varphi \in V_\ell$ satisfies $\varphi(t) > - \frac{1-t}{2}$ for every $t\in (0,T]$, then 
\ben
\label{98}
\lim_n {\cal L}_{n, \lfloor nT\rfloor}^G (\varphi) =\Lambda_{[0,T]}^G (\varphi) := \int_0^T g(t, \varphi(t)) \ dt
\,,
\een
where, for $\theta > - (1-t)/2$
\ben
\label{defg}
g(t, \theta) := \frac{1}{2} \left({\cal J}(1-t +2\theta) -{\cal J}(1-t)- {\cal J}(1 + 2\theta)\right)\,.
\een
\end{lem}
\proof The key point is a convergence of Riemann sums. 
From 
(\ref{LambdaHH}) and (\ref{bin1}) we have, for any $\theta > - \frac{n-k+1}{2n}$,
\ben
\label{ell-ellH}
\Lambda_{n,k}^G (n\theta) &=&
\frac{n-k+2n\theta}{2} \log \left( 1 - \frac{k}{n} + 2 \theta +\frac{1}{n}\right) \\
\nonumber
&-& \frac{n-k}{2} \log \left( 1 - \frac{k}{n} +\frac{1}{n}\right) 
-\frac{n-1 + 2n\theta}{2}\log (1+2\theta)
+ R_{n,k} (\theta)
\een
where
\be
R_{n,k} (\theta) = \int_0^\infty f(s) e^{-\frac{s}{2}}\left[
e^{-\frac{n-k +2n\theta}{2}s} - e^{-\frac{n-k}{2}s}
- e^{-\frac{n-1 +2n\theta}{2}s} + e^{-\frac{n-1}{2}s}\right] ds\,,
\ee
is bounded :
\be
|R_{n,k} (\theta)| \leq 2 \int_0^\infty e^{-\frac{s}{2}} f(s)\!\ ds\,.
\ee
If we set
\be
2\Phi_n (t) &:=& (1-t + 2\varphi(t))\log ( 1 -t + 2\varphi(t) + \frac{1}{n})\\
\nonumber
&& - 
(1 -t) \log( 1-t+ \frac{1}{n})-\Big(1-\frac{1}{n} + 2\varphi(t)\Big)\log (1+2\varphi(t))
\ee
then, making $\theta = \varphi(k/n)$ in (\ref{ell-ellH}), and adding in $k$, we get from (\ref{cgfnG})
\be
\frac{1}{n^2}\Big({\cal L}_{n, \lfloor nt\rfloor}^G (\varphi)- \sum_2^{\lfloor nt\rfloor} R_{n,k}(\varphi)\Big) 
= \frac{1}{n}\sum_1^{\lfloor nt\rfloor}\Phi_n\Big(\frac{k}{n}\Big)
= \int_{1/n}^{\lfloor nt\rfloor/n} \Phi_n \Big(\frac{[nt]}{n}\Big) dt + \frac{1}{n} 
\Phi_n \Big(\frac{\lfloor nt\rfloor}{n}\Big)\,.
\ee
On the one hand, since $\varphi$ is left continuous,  
$\lim_n \Phi_n\Big(\frac{\lfloor nt\rfloor}{n}\Big) =  2 g(t, \varphi(t))$
for every $t\in [0,T]$. On the other hand the following double inequality  holds true:
\be
\nonumber
2\Phi_n (t) &\geq& 
\big(1-t + 2\varphi(t)\big) \log \big(1-t +2 \varphi(t)\big) - (1-t) \log (2-t)
 \\&-& \big(1+ 2\varphi(t)\big)\log \big(1+ 2\varphi(t)\big) -|\log \big(1-t +2 \varphi(t)\big)|
\\
\nonumber
  2 \Phi_n (t)&\leq&
 \big(1-t + 2\varphi(t)\big) \log \big(2-t +2 \varphi(t)\big) - (1-t) \log (1-t) \\
 &-& \big( 1 + 2\varphi(t) \big) \log \big(1+ 2\varphi(t)\big) +|\log \big(1-t +2 \varphi(t)\big)|\,, 
\ee
and with our assumptions on $\varphi$,  these bounds are both integrable. This allows to apply the dominated convergence 
theorem which ends the proof of Lemma
 \ref{cvcgfH}. \QED
If there exists  $s < T$ such that  $2\varphi(s) < - (1-s)$ then for $n$ large enough, 
${\cal L}_{n,\lfloor nT\rfloor}(\varphi) = +\infty$ and we set
$\Lambda_{[0,T]}^G(\varphi) = \infty$. In the other cases we do not know what happens, but
as in \cite{grz}, we will study the exposed points. Before, we need another expression of the dual 
of $\Lambda_{[0,T]}^G$.

\subsubsection{Variational formula}
\label{variaf}

Let us define $\Lambda_{[0,T]}^G (\varphi) = +\infty$ if $\varphi$ does not satisfy the assumption of Lemma \ref{cvcgfH}. 
The dual of $\Lambda_{[0,T]}^G$ is then
\ben
\label{dualeH1}
\left(\Lambda_{[0,T]}^G\right)^\star (\nu) = \sup_{\varphi \in V_\ell} \left\{<\nu, \varphi> - \Lambda_{[0,T]}^G (\varphi)
\right\}
\een
for $\nu \in V_\ell^*$.
Mimicking the method of L\'eonard (\cite{leo1} p. 112-113), we get 
\ben
\label{dualeH2}
\left(\Lambda_{[0,T]}^G\right)^\star (\nu) = \sup_{\varphi \in {\cal C}} \left\{<\nu, \varphi> - \Lambda_{[0,T]}^G (\varphi)
\right\}
\een
where ${\cal C}$ is the set of continuous functions from $[0,T]$ into $`R$ vanishing at $0$.
Then we apply Theorem 5 of Rockafellar \cite{Rocky1}. We get
\be
\left(\Lambda_{[0,T]}^G\right)^\star (\nu) = \int_0^T g^\star\left(t, \frac{d\nu_a}{dt}\right) \!\ dt + \int_0^T r\left(t, \frac{d\nu_s}{d\mu}(t)\right)\!\ d\mu(t)
\ee
where 
\ben
\label{identifH}
g^\star(t,y) &=& \sup_\lambda \left\{\lambda y - g(t,\lambda)\delta(\lambda | (-1/2, \infty)) \right\}\,.
\een
and $r$ is the 
 recession function :
 $$r(t,y) = \lim_{\kappa \rightarrow \infty} \frac{g^\star (t,\kappa y)}{\kappa}\,.$$
 Actually, if $y < 0$, the supremum is achieved for
\ben
\label{lambdaH}
\lambda(t,y) := -\frac{1}{2}\left(1- \frac{t}{1- e^{y}}\right)\,
\een
and we have
 \ben
 \nonumber
g^\star(t,y)&=& \lambda(t,y) y -g \left(t, \lambda(t,y)\right)\\
\nonumber
&=& \frac{1}{2}\Big[-y(1-t) +(1-t) \log(1-t) + t \log t -t \log (1-e^{y})\Big]\\
\label{gstar}
&=& \frac{1}{2} {\bf H}\left(1-t | e^{y}\right)
\,.
\een 
If $y \geq 0$, $g^\star(t,y) 
= \infty$.  The recession is now
 $ r(t,y )= -\frac{1}{2} (1-t) y$ if $y \leq 0$, and $=\infty$ si $y >0$.
As a result
\ben
\label{gLa}
g^\star (t,y) = L_a^G (t,y)\ \ , \ \ r(t,y) = L_s^G (t,y)\,.
\een
 So we proved the identification $\left(\Lambda_{[0,T]}^G\right)^\star = \widetilde I^G_{[0,T]}$
(recall (\ref{tauxHT})).
\subsubsection{Exponential tightness}
\label{expotigh}
If $V_\ell^*$ is provided with the topology $\sigma(V_\ell^*, V_\ell)$, the set 
$B_a := \{\mu \in V_\ell^* : |\mu|_{[0,T]} \leq a\}$ is compact according to the Banach-Alaoglu theorem.
But $\frac{1}{n}\dot\Upsilon_n $ is a positive measure and 
$\frac{1}{n}\dot\Upsilon_n([0,T]) = \frac{1}{n} \Upsilon_n (T) $
has a n.c.g.f. given for  $\theta > 0$ by
\be
\widehat{\cal L}_{n,T}(\theta) := \frac{1}{n^2}\log `E \exp \{n\theta\Upsilon_n (T)
\} 
=  {\cal L}_{n, \lfloor nT\rfloor}(\theta\BBone_{[0,T]}) \ee
For 
 $\theta > - \frac{1-T}{2}$ let
\ben
\label{defLT}
L_T(\theta) := \int_0^T g(t, \theta)\!\ dt\,.
\een 
Lemma \ref{cvcgfH} says that for fixed
 $\theta > - \frac{1-T}{2}$
\ben
\label{bbone} \lim {\cal L}_{n, \lfloor nT\rfloor}(\theta\BBone_{[0,T]}) = L_T(\theta)
\een
so that
\be
\limsup_n \frac{1}{n^2} \log `P \left(\frac{1}{n}\dot\Upsilon_n \notin B_a\right) \leq  
\limsup_n \frac{1}{n^2} \log `P \left(\Upsilon_n (T)  > na \right) 
\leq - a\theta + L_T (\theta) \,,
\ee
which proves the exponential tightness, letting $a \rightarrow \infty$. 

Let us notice that it was not possible to take  $T=1$.

\subsubsection{Exposed points}
\label{exposed}
Let ${\cal R}$ be the set of functions from $[0,T]$ into $`R$ which are positive, continuous and with bounded variation. 
Let ${\cal F}$ be the set of those $m \in V_\ell^*$ (identified with  ${\cal M}_T$ as in \cite{leo1}) 
which are absolutely continuous and whose 
density $\rho$ is such that $-\rho \in {\cal R}.$ 
Let us prove that such a $m$ is exposed, with  exposing hyperplane 
$f_m(t) = \lambda(t, \rho(t))$ (recall (\ref{lambdaH})). Actually we follow the method of  \cite{grz}. For fixed $t$, 
$g^\star(t, .)$ is strictly convex on $(-\infty, 0)$ so that, if $z\not= \rho(t)$, we have
$$g^\star(t, \rho(t)) - g^\star(t, z) < \lambda(t,\rho (t)) (\rho (t) -z)\,. $$
Let 
 $d\xi = \tilde l(t) dt + \xi^\perp$  the Lebesgue decomposition of some element  $\xi \in {\cal M}_T$ such that 
 $\widetilde I^G_{[0,T]}(\xi) < \infty$. 
Taking $z = \tilde l(t)$ and integrating, we get
\be
\int_0^T g^\star (t, \rho(t) dt - \int_0^T g^\star (t, \tilde l(t)) dt < 
\int_0^T \lambda(t,\rho (t)) \rho (t) dt - \int_0^T \lambda(t,\rho (t)) \tilde l(t) \ dt
\ee
and since $\int_0^T g^\star (t, \tilde l(t)) dt = \int_0^T L_a^G (t, \tilde l(t)) dt \leq  \widetilde I_{[0,T]}^G (\xi)$ this yields
\be
\widetilde I_{[0,T]}^G (m) - \widetilde I_{[0,T]}^G (\xi)  < \int_0^T f_m dm -\int_0^T f_m d\xi\,.
\ee
Now let us prove  that this set of exposed points is rich enough. 
We have the following lemma.
\begin{lem}
\label{exposition}
Let $m \in V_r$ such that $\widetilde I_{[0,T]}^G (m) < \infty$. There exists a sequence of functions $l_n \in {\cal R}$ 
 such that 
\begin{enumerate}
\item $\lim_n l_n(t) dt =  -m$ in $V_\ell^*$ with the $\sigma(V_\ell^*, V_\ell)$ topology 
\item  $\lim_n \widetilde I_{[0,T]}^G (-l_n(t) dt) = \widetilde I_{[0,T]}^G(m)\,.$
\end{enumerate}
\end{lem}
\proof
The method may be found in \cite{grz} and in \cite{gamb}. The only difference is in the topology 
because we want to recover marginals.
We will use  the basic inequality which holds for every   $\epsilon \leq 0$ :
\ben
\label{ineqH}
L_a^G (t, v+ \epsilon) \leq L_a^G (t,v) - \frac{\epsilon}{2}(1-t)
\een
Let $m= m_a + m_s$ such that $\widetilde I_{[0,T]}^G (m) < \infty$.  
From (\ref{42}) and (\ref{LaLs}) it is clear that $-m_a$ and $-m_s$ must be  positive measures.
\medskip
 
\noi{First step}
We assume that $m = -l(t) dt - \eta$ with 
$l \in L^1([0,T]; dt)$ and $\eta$ a singular positive measure.
 One can find a sequence of  non negative continuous functions 
 $h_n$  such that $h_n(t) dt \rightarrow \eta$ 
for the topology $\sigma(V_\ell^*, V_\ell)$. Indeed
every function $\psi \in V_\ell$ may be written as a
difference 
 $\psi_1 - \psi_2$ of two increasing functions. There exists a unique (positive) measure $\alpha_{1}$ 
 such that $\psi_1(t) = \alpha_{1} ([t,T])$ for every $t \in [0,T]$. Moreover, 
 the function $g = \eta([0,\cdot]) \in V_r$ is non decreasing and may be approached by a sequence 
  of continuously derivable and non decreasing functions $(g_n)$ such that
  $g_n \leq g$. Setting $h_n := g'_n$ and 
 $\nu_n = h_n (t) dt$, the dominated convergence theorem gives
\be
<\psi_{_1}, \nu_n> = \int_0^T \nu_n ([0,t])\alpha_{1} (dt)  \rightarrow \int_0^T \eta([0,t]) \alpha_{1} (dt) \,.
\ee   
With the same result for  $\psi_2$ we get
\be
<\psi, \nu_n> &=& \int_0^T \nu_n ([0,t])\alpha_{1} (dt) - \int_0^T \nu_n ([0,t])\alpha_{2} (dt) \\&\rightarrow& 
\int_0^T \eta([0,t]) \alpha_{1} (dt) - \int_0^T \eta([0,t]) \alpha_{2} (dt) 
\,.
\ee
or $\lim_n <\psi, \nu_n> = <\psi, \eta>$. 
On the one hand, the lower semicontinuity of 
  $\widetilde I_{[0,T]}^G$ yields
\be
\liminf_n \widetilde I_{[0,T]}^G \Big(-(l(t) + h_n(t))dt\Big)\geq \widetilde I_{[0,T]}^G (m)\,. 
\ee 
On the other hand, integrating (\ref{ineqH}) yields
\be
\widetilde I_{[0,T]}^G (-(l(t) + h_n(t))dt)&\leq& \int_0^T L_a^G (t, -l(t)) dt  + \frac{1}{2}\int_0^T (1-t)h_n (t) dt\\
\nonumber
&\rightarrow& \int_0^T L_a^G (t, -l(t)) dt  + \frac{1}{2}\int_0^T (1-t) \eta (dt) = \widetilde I_{[0,T]}^G (m)\,.
\ee 
\medskip

\noi{Second step} 
Let us assume that $m = -l(t) dt$ with $l \in L^1([0,T]; dt)$ and for every $n$, let us set $l_n = \max(l, 1/n)$. 
It is clear that 
as $n \rightarrow \infty$, then 
$l_n \downarrow l$. 
On the one hand the lower semicontinuity gives $\liminf_A \widetilde I^G_{[0,T]} (-l_n(t)dt) \geq I_{[0,T]}^G (-l(t)dt)$. 
On the other hand, by integration of inequality (\ref{ineqH}), since $l_n -l \leq 1/n$ 
$$I^G_{[0,T]} (-l_n(t)dt) \leq I^G_{[0,T]} (-l(t)dt) + \frac{1}{2n}\,.$$
It is then possible to reduce the problem to the case of functions  bounded below.
\medskip

\noi{Third step}
Let us assume that $m = -l(t) dt$ with $l \in L^1([0,T]; dt)$ and bounded below by $A >0$.
One can find a sequence of continuous  functions $(h_n)$ with bounded variation such that $h_n \geq A/2$
 for every $n$ and such that  
 $h_n \rightarrow l$ a.e. and in $L^1 ([0,T] , dt)$. 
We have $h_n(t) dt \rightarrow l(t)dt$ in $\sigma(V_\ell^*, V_\ell)$ 
and since $L_a^G (t, \cdot)$ 
is uniformly Lipschitz on $]-\infty, -A/2]$, say with constant $\kappa$, we get
\be
|\widetilde I_{[0,T]}^G (-h_n(t) dt) - \widetilde I_{[0,T]}^G (-l(t)dt)|\leq \kappa \int_0^T |h_n (t) - l(t)| dt \rightarrow 0\,.
\ee
Actually, $h_n \in {\cal R}$ and $\varphi_n (t) := \lambda(t, -h_n(t))$ satisfies the assumption of Lemma \ref{cvcgfH} since
$$1 + 2 \varphi_n (t) -t\geq \frac{t}{1-e^{-A/2}}\,.$$

\subsubsection{End of the proof of Theorem \ref{LDPH}}
\label{end}
$\bullet$ First step: upperbound for compact sets. We use th. 4.5.3 b) in \cite{DZ} and the following lemma.
\begin{lem}
\label{petitlemme}
For every $\delta > 0$ and $m \in V_\ell^*$, there exists $\varphi_\delta$ fulfilling conditions of \ref{cvcgfH} 
and such that
\ben
\label{condZani}
\int_0^T \varphi_\delta dm - \Lambda_T^G (\varphi_\delta) \geq \min \Big[ I_{[0,T]}^G (m) - \delta ,\  \delta^{-1}\Big]\,.
\een 
\end{lem}

$\bullet$ Second step: upperbound for closed sets. We use the exponential tightness.

$\bullet$ Third step: lowerbound for open sets. The method is classical (see \cite{DZ} th. 4.5.20 c)), 
owing to Lemma  \ref{exposition}.
\medskip
 
To prove Lemma \ref{petitlemme}, we start from the definition (\ref{dualeH1}) or (\ref{dualeH2}). 
One can find $\bar \varphi_\delta \in V_\ell$ satisfying (\ref{condZani}). If $\bar \varphi_\delta$ 
does not check assumptions of the lemma we add $\varepsilon > 0$ to $\bar \varphi_\delta$ which allows to check them and satisfy 
 (\ref{condZani}) up to a change of $\delta$. \QED 

\subsection{Proof of Theorem \ref{margH}}

We use the contraction from the LDP for paths.
Since the mapping $m \mapsto m([0,T])$ is continuous from $D$ to $\BBr$, 
the family $\{\Upsilon_{n, \lfloor nT\rfloor}\}_n$ satisfies the   LDP with good rate function given by (\ref{opt1}):
$$I^{G}_T (\xi) = \inf \{ I_{[0,T]}^G (v) \ ; \  v(T) = \xi\}\,.$$
Fixing $\xi$, we can look for optimal $v$. Let $\theta  \in (- (1-T)/2, \infty]$ (playing the role of  a Lagrange multiplier). 

By the duality property (\ref{identifH})
\[ 
g^\star \left(t , \frac{d\dot v_a}{dt}(t)\right)\geq \theta \frac{d\dot v_a}{dt}(t) - g(t, \theta)
\,.
\] 
 Integrating and using (\ref{tauxHT}), (\ref{gLa}) and (\ref{defLT}) we get
\ben
I_{[0,T]}^G (v) \geq \theta \dot v_a ([0,T]) - L_T (\theta)- \frac{1}{2}\int_0^T (1-t) \ d\dot v_s (t)\,, 
\een
For every $v$ such that $v(T) = \xi$ it turns out that
\ben
\label{riche} I_{[0,T]}^G (v) \geq \theta \xi -L_T (\theta) - \frac{1}{2}\int_0^T (1-t + 2 \theta) \ d\dot v_s (t)\geq \theta \xi -L_T (\theta)
\,.
\een
Besides, from (\ref{lambdaH}) the ordinary differential equation
\be
\lambda(t, \phi'(t))&=&\theta\\
\phi(0) &=& 0\,,
\ee
admits for unique solution in ${\cal C}^1([0,T])$
\[
t \mapsto \phi (\theta; t) := {\cal J}(1+2\theta)-{\cal J}(1-t +2\theta)  -t \log (1+2\theta)\,.
\]
The mapping $\theta  \mapsto \phi(\theta; T)$
has a positive derivative and its limit as $\theta \downarrow -(1-T)/2$ is $-T$.  
Moreover, by duality
$$g^\star\left(t , \frac{\partial}{\partial t}\phi(\theta,t)\right) = \theta\frac{\partial}{\partial t}\phi(\theta,t) -g(t, \theta)\,.$$
There are two cases.

$\bullet$ If $\xi > - T$, there exists a unique  $\theta_\xi$ such that $\phi(\theta_\xi, T) = \xi$
(i.e. the relation (\ref{=xi1H}) is satisfied). For $v^\xi := \phi(\theta_\xi \!\ , \cdot)$, we get 
from (\ref{tauxHT}), (\ref{gLa}) and (\ref{defLT}) again
\be 
I_{[0,T]}^G(v^\xi) &=& \theta_\xi \xi - L_T(\theta_\xi)
\ee
so that $v^\xi$ realizes the infimum in (\ref{opt1}). A simple computation 
ends the proof of the first statement of Theorem \ref{margH}.

 Let us notice that at the end point $\xi = -T$, we have
 \be
 \theta_\xi = -\frac{1-T}{2} \  , \ v^\xi (t)= {\cal J}(T) -{\cal J}(T-t)  - t \log T\ , \ (v^\xi)'(t) = \log (1-t/T)\,.
 \ee
 Finally
\be
I_T^{G}(-T) &=& T(1-T) +\frac{1}{2}\left(  F(1) - F(1-T)-F(T) + T^2 \log T\right)\\
\nonumber
&=& \frac{T(1-T)}{4} + \frac{T^2 \log T}{4} - \frac{(1-T)^2 \log (1-T)}{4} + \frac{3}{8}
\,.
\ee

$\bullet$ Let us assume $\xi= -T - \varepsilon$ with $\varepsilon > 0$. 
Plugging $\theta = -\frac{1-T}{2}$ in (\ref{riche}) yields, for every $v$ such that $v(T) = \xi$
\be I_{[0,T]}^G (v) \geq - \frac{1-T}{2}\xi -L_T \left(-\frac{1-T}{2}\right) 
=   \frac{1-T}{2} \varepsilon + I_T^{G}(-T)
\ee
Moreover this lower bound is achieved by the measure $\widetilde v = (v^{-T})'(t) dt - \varepsilon \delta_T(t)$, since
\[\int_0^T L_a^G \left(t, (v^{-T})'(t)\right) dt = I_T^{G}(T)\ \ , 
\ \ \int_0^T  \frac{(1-t)}{2}  \!\ \varepsilon \!\ d\delta_T(t) 
=  \frac{(1-T)}{2}\!\ \varepsilon\,.
\]
That ends the proof of the second statement of Theorem \ref{margH}. \QED

\begin{rem}
\label{T=1}
It is possible to try  a direct method to get (\ref{idetH}), (\ref{mauvH}) 
using  G\"artner-Ellis' theorem (\cite{DZ}, Theorem 2.3.6). 
From (\ref{bbone}) the limiting n.c.g.f. of $\frac{\Upsilon_n(T)}{n}$ is $L_T$ 
which is analytic for  $\theta > - \frac{1-T}{2}$.
When $\theta \downarrow -\frac{1-T}{2}$, we have 
$L_T'(\theta)  \downarrow -T$. 
We meet  a case of so called non steepness. To proceed in that direction we could use the method of
 time dependent change of probability
(see \cite{DZT}). We will not give details here. Nevertheless, this approach allows to get one-sided large deviations in the critical case $T=1$.
Actually we get
\be
\lim_n \frac{1}{n^2} \log `P( \Upsilon_{n,n} \geq nx) = - I_1^G (x)
\ee
for $x \geq -1$. The value $x=-1$  corresponds to the limit of $\Upsilon_{n,n}/n$.    
Notice that the second (right) derivative of $I_1^G$ at this point is zero (or equivalently  
$\lim L^{''}_1 (\theta) = \infty$ as $\theta \downarrow 0$)
, which is consistent with previous results on variance. 
I do not know the rate of convergence to $0$ of   $`P( \Upsilon_{n,n} \leq nx)$ for $x< -1$.
\end{rem}
\section{Proofs of Section \ref{Wm}}
\label{sect8}
\subsection{Proofs of Subsection \ref{ascfW}}
We use the decomposition (\ref{HW}). We need only to notice that 
\be
`E \left(\chi_n^2\right)^s &=& 2^s
\frac{\Gamma \left(s+ \frac{n}{2}\right)}{\Gamma \left( \frac{n}{2}\right)}   
\ee
hence
\ben
\label{newcgf}
\log `E \exp \theta \log \frac{\Vert b_1\Vert^2}{n} 
= \ell\Big(\theta + \frac{n}{2}\Big) - \ell\Big(\frac{n}{2}\Big) - \theta \log\Big(\frac{n}{2}\Big)\,, 
\een
which provides estimates for the expectation and the variance. Differentiating once and taking $\theta=0$, we see that
\be
 `E \log  \frac{\Vert b_1\Vert^2}{n} =  \Big[\Psi\Big(\frac{n}{2}\Big) - \log\Big(\frac{n}{2}\Big)\Big]
= -\frac{1}{n}-\int_0^\infty e^{-s\frac{n}{2}} s f(s)\ ds = -\frac{1}{n} + O\Big(\frac{1}{n^{2}}\Big)
\ee
(see (\ref{bin3}), (\ref{proprf})), 
which gives
\ben
\label{espaux} \sup_{p \leq n}\Big|`E S_{n,p} + \frac{p}{n}\Big| = O\Big(\frac{1}{n}\Big)\,.\een
Besides, differentiating (\ref{newcgf}) twice and taking   $\theta=0$ again, we have 
\be
\hbox{Var} \left(\log \frac{\Vert b_1\Vert^2}{n}\right) = \Psi'\Big(\frac{n}{2}\Big) 
= \frac{2}{n} + O\Big(\frac{1}{n^{2}}\Big)
\ee
(see (\ref{restepsi})),  which yields
\ben
\label{varaux}
\sup_{p \leq n} \Big|\hbox{Var} \!\ S_{n,p} - \frac{2p}{n}\Big| = O\Big(\frac{1}{n}\Big)\,.
\een
From (\ref{espaux}) and $(\ref{varaux})$ it is easy to check (via a fourth moment estimate) that  
 $S_n$
 converges in distribution in $\BBd([0,1])$ to 
 $\left(-t + \sqrt{2}\!\ \widetilde{\bf B}_t, \ t \in [0,1]\right)$, where $\widetilde{\bf B}$ 
 is a Brownian motion independent of the $\sigma$-field generated by 
 $(\Upsilon_n , n \in \BBn)$. Finally 
 the family of processes $\eta_n^W = \eta_n^G + S_n$ converges in distribution towards 
 $\big( Y_t -t + \sqrt{2}\!\ \widetilde{\bf B}_t , \ t \in [0,1)\big)$.
 It is a Gaussian process, whose drift and variance are
 \be \d_W(t) = \d_G (t) -t = \frac{1}{2} \log (1-t)\ 
 \,, \ \ 
 \sigma_W^2(t) = \sigma_G^2 (t) + 2t = -2 \log (1-t)\,.
 \ee
 which identify the process $X$. 
 
 Besides, we have
 $$\widehat\eta^W_n (1) =  \widehat\eta^G_n (1) + \frac{S_n(1)}{\sqrt{2\log n}}\,,$$
 so that the convergence of $\widehat\eta^W_n (1)$ is clear. Moreover the independence properties seen 
 in Section \ref{flucH} remain true.
\subsection{Proofs of Subsection \ref{ldpW}}
\label{similarroute}
\subsubsection{Proof of Lemma \ref{LDPaux}}
It is a route similar to the proof of Theorem \ref{LDPH} in 
Section \ref{7.4} (see also \cite{Najim1}). We start from (\ref{newcgf}):
\be
\log `E  \exp< n \dot S_n , \gamma >
  = \sum_{k=1}^n
  \left[ n \gamma\left(\frac{k}{n}\right) \log 2 + \ell\left(n \gamma\left(\frac{k}{n}\right) + \frac{n}{2}\right) - 
  \ell\left(\frac{n}{2}\right)\right]
 \ee
if $\gamma(s) + \frac{1}{2} > 0$ for every $s \in [0,1]$.
The limiting n.c.g.f. is 
\ben
\label{128}
{\cal L}^{S} (\gamma) = \frac{1}{2} \int_0^1 {\cal J}(1+2\gamma(t)) dt\,,
\een
which yields (\ref{60}) by duality (see \cite{Rocky1} again).
\subsubsection{Proof of Theorem  \ref{LDPwish}}We deduce from Lemma \ref{LDPaux} and Theorem \ref{LDPH} 
that the sum $\frac{1}{n}\dot \Upsilon_n + \frac{1}{n}\dot S_n$ 
satisfies a LDP at the same scale with good rate function obtained by 
inf-convolution of $I_{[0,T]}^{G}$ and $\widetilde I_{[0,T]}^{S}$. To obtain (\ref{a+s}) and (\ref{identif}), 
it is possible to compute explicitely this inf-convolution:
\begin{eqnarray*}
L_a^W
 &=& \inf_v \{L_a ^{G}(v) + L_a^{S}(u-v)\}\\
L_s^W
&=& \inf_v \{L_s^{G}(v) + L_s^{S}(u-v)\}\,.
\end{eqnarray*}
Alternatively, it is possible to sum the two n.c.g.f. ((\ref{98}) and (\ref{128})) and get the rate function by  duality. \QED
\begin{rem}
We can make the same comments as in Remark \ref{T=1}. In particular, we get
\ben
\lim_n \frac{1}{n^2} \log `P( \Xi_{n,n} \geq nx) = - I_1^W (x)
\een
for $x \geq -1$. This boundary point corresponds to the limit of $\Xi_{n,n}/n$.  The second derivative is vanishing at this point,
which is consistent with  the results on variance.
\end{rem}
\subsubsection{Proof of Proposition  \ref{infI}}Let $\theta$  be a Lagrangian factor. We begin  by minimizing
\be
I_T^{sp} (\mu)- \theta T \int \log x \ d\mu(x)  = \frac{T^2}{2} 
\left[- \Sigma(\mu) + 2 \int q_{\lambda,s} (x)
 \ d\mu(x)\right] + B(T)
\ee
where
\ben
\label{lambdas}
q_{\lambda, s} (x) = \lambda x - s \log x, \ \ \lambda = \frac{1}{2T} , \ s = \frac{1-T +2 \theta}{2T}\,.
\een
In \cite{Saff} p.43 example 5.4,
it is stated that  for $\lambda > 0$ and $2s+1 > 0$ fixed, the infimum $$\inf_\mu \ - \Sigma(\mu) + 2 \int q_{\lambda,s} (x)
\ d\mu(x)$$ 
  is achieved on the unique
extremal measure 
$\pi_{\sigma^2}^c$ 
 with 
\be
\sigma^2 = \frac{2s+1}{2\lambda} \ , \ c = \frac{1
}{2s+1}\,,\ee
which yields, from (\ref{lambdas}):
\be
\sigma^2 = 1 +2\theta , \ \ \ c = \frac{T}{\sigma^2}\,.
\ee
Now it remains to look for $\theta$ such that the constraint $\mu \in {\cal A}(\xi/T)$ is saturated.
Since
\be
\int \log x \ d\pi_{\sigma^2}^c (x) = \log \sigma^2 + \int \log x \ d\pi_1^c(x) dx\,, 
\ee
and thanks to (\ref{often}) we see that $\theta$ must satisfy 
\be
\xi = T \log \sigma^2 - T \frac{{\cal J}(1-c)}{c} = {\cal J}(1+2\theta) - {\cal J}(1-T + 2\theta)\,,
\ee
which is exactly exactly (\ref{=xi1W}). 

To compute $I_T^{(sp)}(\pi_{\sigma^2}^c)$, we start from the definition (\ref{HP1}):
\be
I_T^{(sp)}(\pi_{\sigma^2}^c) = -\frac{T^2}{2}\Sigma(\pi_{\sigma^2}^c) + \frac{T}{2}\int (x - (1-T)\log x)\ 
d\pi_{\sigma^2}^c(x) + B(T)\,,
\ee
and transform $\pi_{\sigma^2}^c$ to $\pi_1^c$ using the dilatation. In particular,  
(\ref{HP2}) yields $\Sigma(\pi_{\sigma^2}^c) = \log \sigma^2 + \Sigma(\pi_1^c)$ and $\Sigma(\pi_1^c)$ 
may be picked from formula (13) p.10 in \cite{hiai1} :
$$\Sigma(\pi_1^c) = -1 + \frac{1}{2}\left(c^{-1} + \log c\ + (c^{-1} -1)^2 \log(1-c)\right)\,.$$ 
Besides   we have easily $\int x \ d\pi_1^c (x) = 1$.
After some tedious but elementary computations we get exactly the expression (\ref{idetW}).
\QED
\section{Appendix : Some properties of $\ell = \log \Gamma$ and $\Psi$}
\label{appendix}
From the Binet's formula (\cite{astig} \cite{bateman} p.21), we have
\ben
\label{bin1}
\ell (x) &=& (x-\frac{1}{2})\log x -x +1 + \int_0^\infty f(s)[e^{-sx} - e^{-s}]\ \! ds\\
\label{bin2}
&=& (x-\frac{1}{2})\log x -x + \frac{1}{2} \log(2\pi) + \int_0^\infty f(s)e^{-sx}\ \! ds\,.
\een
where the function $f$ is defined by
\ben
\label{propf}
 f(s) = \left[\frac{1}{2}-\frac{1}{s}+ \frac{1}{e^s -1}\right]\frac{1}{s} = 
 2\sum_{k=1}^\infty \frac{1}{s^2 + 4\pi^2 k^2}\,,
\een
and satisfies for every $s \geq 0$
\ben
\label{proprf}
0 < f(s) \leq f(0)= 1/12 \ , \  \ 0 < \big(sf(s)+ \frac{1}{2}\big) < 1\,.
\een
By differentiation
\ben
\label{bin3}
\Psi(x) = \log x - \frac{1}{2x} -\int_0 ^\infty s f(s) e^{-sx}\, ds = 
 \log x -\int_0 ^\infty e^{-sx} \big(sf(s)+ \frac{1}{2}\big) ds\,, 
\een
As easy consequences, we have, for every $x > 0$
\ben
\label{supx}
0 < x\left(\log x - \Psi (x)\right) \leq 1\,, 
\een
 and
\ben
\label{supx2}
0 < x^2\left(\log x - \Psi(x)- \frac{1}{2x}\right) \leq \frac{1}{12}\,.
\een
Differentiating again we see that
for $q\geq 1$ 
\ben
\label{polygamma}
\Psi^{(q)}(z) = (-1)^{q-1} q! z ^{-q} + (-1)^{q-1} \int_0^\infty e^{-sz} s^q  \big(sf(s)+ \frac{1}{2}\big)\ ds
\een
and then
\ben
\label{restepsi}
 |\Psi^{(q)}(z) - (-1)^{q-1} q! z ^{-q}| \leq z^{-q-1} q!\,.
 \een

\bibliographystyle{plain}
\bibliography{matnew7}
\end{document}